\documentclass[review,hidelinks,onefignum,onetabnum]{siamart220329}
\usepackage{amsmath}
\usepackage{lipsum}
\usepackage{amsfonts}
\usepackage{graphicx}
\usepackage{epstopdf}
\usepackage{algorithmic}
\usepackage{amssymb}
\usepackage{amsopn}

\nolinenumbers

\ifpdf
\DeclareGraphicsExtensions{.eps,.pdf,.png,.jpg}
\else
\DeclareGraphicsExtensions{.eps}
\fi

\newsiamremark{remark}{Remark}
\newsiamremark{hypothesis}{Hypothesis}
\crefname{hypothesis}{Hypothesis}{Hypotheses}
\newsiamthm{claim}{Claim}
\headers{Convexification for Mean Field Games}
{M. V. Klibanov, J. Li, and Z. Yang}

\ifpdf
\hypersetup{
  pdftitle={Convexification for Mean Field Games},
  pdfauthor={M. V. Klibanov, J. Li, and Z. Yang}
}
\fi
\input{tcilatex}
\begin{document}

\title{Convexification Numerical Method for the Retrospective Problem of
Mean Field Games \thanks{%
Submitted Date. 
\funding{The work of J. Li was partially
supported by the NSF of China No. 11971221, Guangdong NSF Major Fund No.
2021ZDZX1001, the Shenzhen Sci-Tech Fund No. RCJC20200714114556020,
JCYJ20200109115422828 and JCYJ20190809150413261.} }}
\author{Michael V. Klibanov
	\thanks{Department of Mathematics and Statistics, University of North Carolina at Charlotte, Charlotte, NC, 28223, USA 
		(\email{mklibanv@charlotte.edu}).}
	\and Jingzhi Li
	\thanks{Department of Mathematics \& National Center for Applied Mathematics Shenzhen \& SUSTech International Center for Mathematics, Southern University of Science and Technology, Shenzhen 518055, P.~R.~China 
		(\email{li.jz@sustech.edu.cn}).}
	\and Zhipeng Yang
	\thanks{Department of Mathematics, Southern University of Science and Technology, Shenzhen 518055, P.~R.~China 
		(\email{yangzp@sustech.edu.cn}).}
}
\maketitle

\begin{abstract}
The convexification numerical method with the rigorously established global
convergence property is constructed for a problem for the Mean Field Games
System of the second order. This is the problem of the retrospective
analysis of a game of infinitely many rational players. In addition to
traditional initial and terminal conditions, one extra terminal condition is
assumed to be known. Carleman estimates and a Carleman Weight Function play
the key role. Numerical experiments demonstrate a good performance for
complicated functions. Various versions of the convexification have been
actively used by this research team for a number of years to numerically
solve coefficient inverse problems.
\end{abstract}




\section{Introduction}

\label{sec:1}

The mean field games (MFG) theory studies the collective behavior of
infinitely many rational agents. This theory was introduced in the seminal
publications of Lasry and Lions \cite{LL} as well as of Huang, Caines, and
Malham\'{e} \cite{Huang}. Social sciences play an increasingly significant
role in the modern society. Therefore, mathematical modeling of social
phenomena has the potential for a substantial societal impact. The MFG
theory stands out as the only mathematical model of social processes that
relies on a universal system of coupled Partial Differential Equations
(PDEs) known as the Mean Field Games system (MFGS). The MFGS is the key
ingredient of the MFG theory \cite{A}. The MFGS of the second order is a
system of two coupled nonlinear parabolic PDEs with two opposite directions
of time. The applications of this theory to societal problems are broad and
encompass such areas as, e.g. finance, combating corruption, cybersecurity,
interactions of electrical vehicles, quantum information theory, election
dynamics, robotic control, etc., see, e.g. \cite{A,Co,Hui,KM,Kol,Trusov} for
a far incomplete list of references. Another notable application is in
boundary control problems \cite{Carmona2018}.

Hence, due to a wide range of applications of the MFG theory, it becomes
important to address various mathematical questions regarding the MFGS. In
this paper we present the \underline{first} numerical method for the MFGS
with the rigorously guaranteed global convergence property. We address the
application to the retrospective analysis of games, which are governed by
MFGS. In other words, our numerical method is aimed to figure out the
history of a mean field game after that game is finished. 

This paper consists of two equally important parts. The first part is
theoretical and the second part is numerical. Numerical studies of the
second part confirm the theory of this first one.

\textquotedblleft Global convergence" means that our convergence analysis
ends up with a theorem, which guarantees that the iterative solutions,
generated by our method, converge to the true solution of the MFGS (if it
exists) starting from any point of an \emph{a priori} given convex bounded
set in a Hilbert space. The diameter $d>0$ of this set is, though fixed, but
an arbitrary number. That convergence takes place as long as the level of
noise $\delta >0$ in the input data tends to zero. An explicit estimate of
the convergence rate is given. Results of our numerical experiments
demonstrate a good accuracy of computed solutions in the presence of the
random noise in the input data. We refer to, e.g. \cite{Co,LiuOsher,Trusov}
and references cited therein for some other numerical methods for the MFGS.

In this paper we develop a version of the convexification method for the
problem of the retrospective analysis of the processes governed by the MFGS
of the second order. In the past the convexification was applied by this
research group for constructions and numerical testing of globally
convergent numerical methods for some ill-posed Cauchy problems for
nonlinear PDEs as well as for some coefficient inverse problems for PDEs.
The latter problems are also nonlinear. See, e.g. \cite{Bak,LeLoc1,LeLoc2}, 
\cite[chapter 5]{KL} for selected samples of publications on the first named
topic and \cite[Chapters 7-12]{KL}, \cite{SAR,Kpar,Ktransp,LeLoc2} for the
second topic; also, see references cited therein.

Our problem was first considered in the work of Klibanov and Averboukh \cite%
{MFG1}. The work \cite{MFG1} is the first one, in which the tool of Carleman
estimates was introduced in the MFG theory. The Lipschitz stability estimate
was proven in \cite{MFG1}. In the follow up publications Carleman estimates
were applied to prove both H\"{o}lder and Lipschitz stability estimates for
various problems for the MFGS, see \cite{ImLY,MFG2,MFG6} and references
cited in \cite{MFG6}. These stability estimates can also be interpreted as
the accuracy estimates of the solutions of the MFGS with respect to the
noise in the input data. Those stability estimates imply uniqueness of
corresponding problems. Stability estimates for the MFGS were unknown prior
these works.

In addition, the framework of the regularization theory for Ill-Posed and
Inverse Problems was introduced in the MFG theory in \cite{MFG1}. This
framework is also used in \cite{ImLY,MFG2,MFG6} as well as in the current
paper.

\textbf{Remark 1.1.} \emph{We now outline one of fundamental principles of
the regularization theory for Ill-Posed and Inverse Problems \cite{T}, which
is used in the publications of this group on the MFG theory. It is assumed
that the input data for a problem under consideration are noisy with the
level of noise }$\delta >0$\emph{. It is further assumed that there exists
an \textquotedblleft ideal" exact/true solution for the \textquotedblleft
ideal" noiseless data. Then the regularization method means that a }$\delta
- $\emph{dependent family of approximate solutions is constructed and these
solutions converge to the true solution as long as }$\delta \rightarrow 0.$

\textbf{Remark 1.2. }\emph{Traditionally minimal smoothness requirements are
not of a significant concern in the theory of Ill-Posed and Inverse
Problems, see, e.g. \cite{Nov}, \cite[Theorem 4.1]{Rom}. Therefore, we are
also not concerned with such requirements here.}

All functions considered below are real valued ones. In section 2 we state
the problem. In section 3 we formulate two Carleman estimates. In section 4
we present the convexification method and carry out its convergence
analysis. Numerical experiments are described in section 5. Summary is
provided in section 6.

\section{Statement of the Problem}

\label{sec:2}

Denote $x=\left( x_{1},...,x_{n}\right) \in \mathbb{R}^{n}$ points in $%
\mathbb{R}^{n}.$ Let $\Omega \subset \mathbb{R}^{n}$ be a bounded domain
with the piecewise smooth boundary $\partial \Omega .$ Let $T>0$ be a
number. Denote $Q_{T}=\Omega \times \left( 0,T\right) ,$ $S_{T}=\partial
\Omega \times \left( 0,T\right) ,$ $\partial _{\nu }v$ normal derivative of
an appropriate function $v\left( x,t\right) $ at $S_{T}.$ Let $s\in
C^{1}\left( \overline{Q}_{T}\right) $ and $K\left( x,y\right) \in L_{\infty
}\left( \Omega \times \Omega \right) $ be some functions, and $\beta >0$ be
a number. The MFGS is \cite{A,LL}:

\begin{equation}
\left. 
\begin{array}{c}
L_{1}\left( u,p\right) =u_{t}(x,t)+\beta \Delta u(x,t)+s(x,t)(\nabla
u(x,t))^{2}/2+ \\ 
+\int\limits_{\Omega }K\left( x,y\right) p\left( y,t\right) dy+f\left(
x,t\right) p\left( x,t\right) +F_{1}\left( x,t\right) =0,\text{ }\left(
x,t\right) \in Q_{T}, \\ 
L_{2}\left( u,p\right) =p_{t}(x,t)-\beta \Delta p(x,t)+\func{div}%
(s(x,t)p(x,t)\nabla u(x,t))+ \\ 
+F_{2}\left( x,t\right) =0,\text{ }\left( x,t\right) \in Q_{T},%
\end{array}%
\right.  \label{2.1}
\end{equation}%
where $L_{1}\left( u,p\right) $ and $L_{2}\left( u,p\right) $ are two
operators. Just as in \cite{MFG1}, we add the zero Neumann boundary
conditions%
\begin{equation}
\partial _{\nu }u\mid _{S_{T}}=\partial _{\nu }p\mid _{S_{T}}=0.  \label{2.2}
\end{equation}%
Conventional initial and terminal conditions for the MFGS (\ref{2.1}) are:%
\begin{equation}
u\left( x,T\right) =u_{T}\left( x\right) ,\text{ }p\left( x,0\right)
=p_{0}\left( x\right) ,\text{ }x\in \Omega .  \label{2.3}
\end{equation}

In (\ref{2.1}) $x$ is the position of an agent, $u\left( x,t\right) $ is the
value function and $p\left( x,t\right) $ is the distribution of agents at
the point $x$ and at the moment of time $t$. The integral operator in (\ref%
{2.1}) is the so-called \textquotedblleft global interaction term", $f\left(
x,t\right) p\left( x,t\right) $ is the local interaction term, and $\beta
\Delta u$ and $\beta \Delta p$ are viscosity terms. We have added functions $%
F_{1},F_{2}$ in (\ref{2.1}) for two reasons. First, to figure out the
accuracy of the minimizer of our functional, see Theorem 4.2 in section 4.
Second, we need $F_{1}$ for the numerical data generation in section 5. Let $%
M>0$ be a number. We assume that%
\begin{equation}
\left. 
\begin{array}{c}
s\left( x,t\right) \in C^{1}\left( \overline{Q}_{T}\right) ;\text{ }f\left(
x,t\right) ,F_{1}\left( x,t\right) ,F_{1}\left( x,t\right) \in L_{\infty
}\left( Q_{T}\right) , \\ 
\left\Vert s\right\Vert _{C^{1}\left( \overline{Q}_{T}\right) },\left\Vert
f\right\Vert _{L_{\infty }\left( Q_{T}\right) },\left\Vert F_{1}\right\Vert
_{L_{\infty }\left( Q_{T}\right) },\left\Vert F_{2}\right\Vert _{L_{\infty
}\left( Q_{T}\right) },\left\Vert K\right\Vert _{L_{\infty }\left( \Omega
\times \Omega \right) }\leq M.%
\end{array}%
\right.  \label{2.4}
\end{equation}

The goal of this paper is to solve numerically the following problem:

\textbf{Problem}. \emph{Find the pair of functions }$u,p\in H^{2}\left(
Q_{T}\right) $\emph{\ satisfying conditions (\ref{2.1})-(\ref{2.3}) as well
as the following terminal condition:}%
\begin{equation}
p\left( x,T\right) =p_{T}\left( x\right) ,\text{ }x\in \Omega .  \label{2.5}
\end{equation}

This problem can be considered as a problem of the retrospective analysis.
In other words, a game is finished at the moment of time $\left\{
t=T\right\} .$ We measure the final distribution of agents $p\left(
x,T\right) .$ And our goal is to figure out the history of this game, i.e.
to get a knowledge on how the game proceeded on the time interval $t\in
\left( 0,T\right) .$

It is well known that uniqueness of the solution of the system (\ref{2.1})
cannot be guaranteed if only conditions (\ref{2.2}) and (\ref{2.3}) are
known, unless strong assumptions are imposed \cite{LL}. On the other hand,
adding condition (\ref{2.5}) guarantees Lipschitz stability estimate and,
therefore, uniqueness of problem (\ref{2.1})-(\ref{2.3}), (\ref{2.5}) \cite%
{MFG1}.

\section{Carleman Estimates}

\label{sec:3}

Carleman estimates play the key role in our technique. Carleman estimates
are independent on low order terms of PDE operators \cite[Lemma 2.1.1]{KL}
and are, therefore, formulated only for principal parts of those operators.
Since two different directions of time are present in two parabolic
equations (\ref{2.1}), then we formulate in this section two Carleman
estimates: for operators $\partial _{t}+\beta \Delta $ and $\partial
_{t}-\beta \Delta .$

Let $a>0$ be a parameter. First, consider the function $\psi _{\mu ,\lambda
}\left( t\right) =e^{2\mu \left( t+a\right) ^{\lambda }},t\in \left(
0,T\right) ,\mu ,\lambda >0.$ This function was used in Carleman estimates
of analytical studies in \cite{MFG1}. However, this function is inconvenient
to work with in computations because it depends on two parameters $\mu
,\lambda .$ Fortunately, Carleman estimates of Theorems 3.1 and 3.2 of \cite%
{MFG1} are formulated and proved for all values $\mu >0.$ Therefore, we use
below the Carleman Weight Function $\varphi _{\lambda }\left( t\right) =\psi
_{1,\lambda }\left( t\right) ,$ i.e. 
\begin{equation}
\varphi _{\lambda }\left( t\right) =e^{2\left( t+a\right) ^{\lambda }},\text{
}t\in \left( 0,T\right) ,\text{ }\lambda \geq 1.  \label{3.1}
\end{equation}

\textbf{Remark 3.1.}\emph{\ At the time of the submission of this
manuscript, the work \cite{MFG1} from which Theorems 3.1 and 3.2 are taken,
is posted online }www.arxiv.org\emph{, although not yet published in a
journal. However, similar Carleman estimates are fully proven in the paper 
\cite{MFG2}, which is currently published. The only difference is in the
Carleman Weight Function, which is chosen in \cite{MFG2} as }$\rho _{\lambda
}\left( t\right) =e^{2\left( T-t+a\right) ^{\lambda }}$\emph{\ instead of
the function }$\varphi _{\lambda }\left( t\right) $\emph{\ in (\ref{3.1}).
This difference is due to the fact that the problem considered in \cite{MFG2}
is different from the problem considered both in \cite{MFG1} and here.
Proofs of Carleman estimates in \cite{MFG2} are quite similar ones with the
proofs of Carleman estimates in \cite{MFG1}.}

Denote 
\begin{equation*}
\left. 
\begin{array}{c}
H_{0}^{2}\left( Q_{T}\right) =\left\{ u\in H^{2}\left( Q_{T}\right)
:\partial _{\nu }u\mid _{S_{T}}=0\right\} , \\ 
\hspace{-0.1cm}H^{1,1}(Q_{T})=\left\{ u:\left\Vert u\right\Vert
_{H^{1,1}\left( Q_{T}\right) }^{2}=\left\Vert u_{t}\right\Vert _{L_{2}\left(
Q_{T}\right) }^{2}+\sum\limits_{i=1}^{n}\left\Vert u_{x_{i}}\right\Vert
_{L_{2}\left( Q_{T}\right) }^{2}+\left\Vert u\right\Vert _{L_{2}\left(
Q_{T}\right) }^{2}<\infty \right\} , \\ 
H^{1,0}\left( Q_{T}\right) =\left\{ u:\left\Vert u\right\Vert
_{H^{1,0}\left( Q_{T}\right) }^{2}=\sum\limits_{i=1}^{n}\left\Vert
u_{x_{i}}\right\Vert _{L_{2}\left( Q_{T}\right) }^{2}+\left\Vert
u\right\Vert _{L_{2}\left( Q_{T}\right) }^{2}<\infty \right\} .%
\end{array}%
\right.
\end{equation*}

\textbf{Theorem 3.1 }(the first Carleman estimate \cite{MFG1})\textbf{.} 
\emph{There exists a number }$C=C\left( T,\beta ,a\right) >0$\emph{\
depending only on listed parameters such that the following Carleman
estimate holds:}%
\begin{equation}
\left. 
\begin{array}{c}
\int\limits_{Q_{T}}\left( u_{t}+\beta \Delta u\right) ^{2}\varphi _{\lambda
}^{2}dxdt\geq C\int\limits_{Q_{T}}\left( u_{t}^{2}+\left( \Delta u\right)
^{2}\right) \varphi _{\lambda }^{2}dxdt+ \\ 
+C\lambda \int\limits_{Q_{T}}\left( \nabla u\right) ^{2}\varphi _{\lambda
}^{2}dxdt+C\lambda ^{2}\int\limits_{Q_{T}}u^{2}\varphi _{\lambda }^{2}dxdt,
\\ 
\forall \lambda >2,\forall u\in H_{0}^{2}\left( Q_{T}\right) \cap \left\{
u\left( x,T\right) =0\right\} .%
\end{array}%
\right.  \label{3.2}
\end{equation}

\textbf{Theorem 3.2 }(the second Carleman estimate \cite{MFG1}). \emph{There
exist a sufficiently large number }$\lambda _{0}=\lambda _{0}\left( T,\beta
,a\right) >2$\emph{\ and a number }$C=C\left( T,\beta ,a\right) >0$\emph{\
depending only on listed parameters such that the following Carleman
estimate holds: } 
\begin{equation}
\left. 
\begin{array}{c}
\int\limits_{Q_{T}}\left( u_{t}-\beta \Delta u\right) ^{2}\varphi _{\lambda
}dxdt\geq C\sqrt{\lambda }\int\limits_{Q_{T}}\left( \nabla u\right)
^{2}\varphi _{\lambda }dxdt+C\lambda ^{2}\int\limits_{Q_{T}}u^{2}\varphi
_{\lambda }dxdt, \\ 
\forall \lambda \geq \lambda _{0},\text{ }\forall u\in H_{0}^{2}\left(
Q_{T}\right) \cap \left\{ u\left( x,0\right) =u\left( x,T\right) =0\right\} .%
\end{array}%
\right.  \label{3.3}
\end{equation}

\section{Convexification}

\label{sec:4}

\subsection{The minimization problem}

\label{sec:4.1}

Define the integer $k_{n},$ 
\begin{equation}
k_{n}=\left[ \left( n+1\right) /2\right] +3,  \label{4.01}
\end{equation}%
where $\left[ \left( n+1\right) /2\right] $ is the largest integer not
exceeding $\left( n+1\right) /2,$ see Remark 1.2. By embedding theorem and (%
\ref{4.01}) $H^{k_{n}}\left( Q_{T}\right) \subset C^{2}\left( \overline{Q}%
_{T}\right) ,$ and there exists a constant $C_{0}=C_{0}\left( Q_{T}\right)
>0 $ depending only on the domain $Q_{T}$ such that%
\begin{equation}
\left\Vert g\right\Vert _{C^{2}\left( \overline{Q}_{T}\right) }\leq
C_{0}\left\Vert g\right\Vert _{H^{k_{n}}\left( Q_{T}\right) },\text{ }%
\forall g\in H^{k_{n}}\left( Q_{T}\right) .  \label{4.3}
\end{equation}

Let $R>0$ be an arbitrary number. Consider the set $B(R),$%
\begin{equation}
B\left( R\right) =\left\{ 
\begin{array}{c}
\left( u,p\right) \in H^{k_{n}}\left( Q_{T}\right) \times H^{k_{n}}\left(
Q_{T}\right) :u,p\in H_{0}^{2}\left( Q_{T}\right) , \\ 
u\left( x,T\right) =u_{T}\left( x\right) ,\text{ }p\left( x,T\right)
=p_{T}\left( x\right) ,\text{ }p\left( x,0\right) =p_{0}\left( x\right) , \\ 
\left\Vert u\right\Vert _{H^{k_{n}}\left( Q_{T}\right) },\left\Vert
p\right\Vert _{H^{k_{n}}\left( Q_{T}\right) }<R.%
\end{array}%
\right\}  \label{4.4}
\end{equation}%
Let $L_{1}\left( u,p\right) $ and $L_{2}\left( u,p\right) $ be two operators
defined in (\ref{2.1}). Consider four functionals 
\begin{equation}
\left. 
\begin{array}{c}
J_{1,\lambda },J_{2,\lambda },J_{3},J:H^{k_{n}}\left( Q_{T}\right) \times
H^{k_{n}}\left( Q_{T}\right) \rightarrow \mathbb{R}, \\ 
J_{1,\lambda }\left( u,p\right) =\int\limits_{Q_{T}}\left( L_{1}\left(
u,p\right) \right) ^{2}\varphi _{\lambda }dxdt,\text{ }J_{2,\lambda }\left(
u,p\right) =\int\limits_{Q_{T}}\left( L_{2}\left( u,p\right) \right)
^{2}\varphi _{\lambda }dxdt, \\ 
J_{3}\left( u,p\right) =\gamma \left( \left\Vert u\right\Vert
_{H^{k_{n}}\left( Q_{T}\right) }^{2}+\left\Vert p\right\Vert
_{H^{k_{n}}\left( Q_{T}\right) }^{2}\right) , \\ 
J_{\lambda ,\gamma }\left( u,p\right) =J_{1,\lambda }\left( u,p\right)
+\left( 1/2+C_{1}/\lambda ^{2}\right) J_{2,\lambda }\left( u,p\right)
+J_{3}\left( u,p\right) ,%
\end{array}%
\right.  \label{4.5}
\end{equation}%
where operators $L_{1}\left( u,p\right) $ and $L_{2}\left( u,p\right) $ are
defined in (\ref{2.1}) and $\gamma \in \left( 0,1\right) $ is the
regularization parameter. To solve our target Problem, we consider below

\textbf{Minimization Problem}. \emph{Minimize the functional }$J_{\lambda
,\gamma }\left( u,p\right) $\emph{\ in (\ref{4.5}) on the set }$B\left(
R\right) $\emph{\ defined in (\ref{4.4}).}

\subsection{The strong convexity of the functional $J_{\protect\lambda , 
\protect\alpha }\left( u,p\right) $ on $B\left( R\right) $}

\label{sec:4.2}

Below $\left[ ,\right] $ is the scalar product in the Hilbert space $%
H^{k_{n}}\left( Q_{T}\right) \times H^{k_{n}}\left( Q_{T}\right) .$ Define
the subspace $\widetilde{H}$ of this space as 
\begin{equation}
\left. \widetilde{H}=\left\{ 
\begin{array}{c}
\left( h,q\right) \in H^{k_{n}}\left( Q_{T}\right) \times H^{k_{n}}\left(
Q_{T}\right) :h,q\in H_{0}^{2}\left( Q_{T}\right) , \\ 
h\left( x,T\right) =\text{ }q\left( x,T\right) =q\left( x,0\right) =0, \\ 
\left\Vert \left( h,q\right) \right\Vert _{\widetilde{H}}^{2}=\left\Vert
h\right\Vert _{H^{k_{n}}\left( Q_{T}\right) }^{2}+\left\Vert q\right\Vert
_{H^{k_{n}}\left( Q_{T}\right) }^{2}.%
\end{array}%
\right\} \right.  \label{4.12}
\end{equation}

\textbf{Theorem 4.1} (the central analytical result).\emph{\ Assume that
conditions (\ref{2.4}), (\ref{4.01}) and (\ref{4.5}) hold. Let }$R>0$\emph{\
be an arbitrary number and }$B\left( R\right) $\emph{\ be the set defined in
(\ref{4.4}). Then:}

\emph{1. The functional }$J_{\lambda ,\gamma }$\emph{\ has the Fr\'{e}chet
derivative }$J_{\lambda ,\gamma }^{\prime }\left( u,p\right) \in \widetilde{H%
}$\emph{\ at every point }$\left( u,p\right) \in \overline{B\left( R\right) }
$\emph{. The Fr\'{e}chet derivative }$J_{\lambda ,\gamma }^{\prime }\left(
u,p\right) $ \emph{is Lipschitz continuos on }$\overline{B\left( R\right) },$%
\emph{\ i.e. the following inequality holds:}%
\begin{equation}
\left. 
\begin{array}{c}
\left\Vert J_{\lambda ,\gamma }\left( u_{1},p_{1}\right) -J_{\lambda ,\gamma
}\left( u_{2},p_{2}\right) \right\Vert _{H^{k_{n}}\left( Q_{T}\right) \times
H^{k_{n}}\left( Q_{T}\right) }\leq \\ 
\leq D\left\Vert \left( u_{1},p_{1}\right) -\left( u_{2},p_{2}\right)
\right\Vert _{H^{k_{n}}\left( Q_{T}\right) \times H^{k_{n}}\left(
Q_{T}\right) },\text{ }\forall \left( u_{1},p_{1}\right) ,\left(
u_{2},p_{2}\right) \in \overline{B\left( R\right) },%
\end{array}
\right.  \label{4.6}
\end{equation}%
\emph{where the number }$D=D\left( \lambda ,\gamma ,\Omega ,T,M,R\right) >0$%
\emph{\ depends only on listed parameters.}

\emph{2. There exists a sufficiently large number }$\overline{\lambda }=%
\overline{\lambda }\left( \Omega ,T,M,R\right) >1$\emph{\ such that for all }%
$\lambda \geq \overline{\lambda }$\emph{\ the functional }$J_{\lambda
,\gamma }$\emph{\ is strongly convex on the set }$\overline{B\left( R\right) 
},$\emph{\ i.e. there exists a number }$C_{1}=C_{1}\left( \Omega
,T,M,R\right) >0$\emph{\ such that the following inequality holds:}%
\begin{equation}
\left. 
\begin{array}{c}
J_{\lambda ,\gamma }\left( u_{1},p_{1}\right) -J_{\lambda ,\gamma }\left(
u,p\right) -\left[ J_{\lambda ,\gamma }^{\prime }\left( u,p\right) ,\left(
u_{1}-u,p_{1}-p\right) \right] \geq \\ 
\geq C_{1}e^{2a^{\lambda }}\left( \left\Vert \Delta u_{1}-\Delta
u\right\Vert _{L_{2}\left( Q_{T}\right) }^{2}+\left\Vert u_{1}-u\right\Vert
_{H^{1,1}\left( Q_{T}\right) }^{2}+\left\Vert p_{1}-p\right\Vert
_{H^{1,0}\left( Q_{T}\right) }^{2}\right) + \\ 
+\gamma \left( \left\Vert u_{1}-u\right\Vert _{H^{k_{n}}\left( Q_{T}\right)
}^{2}+\left\Vert p_{1}-p\right\Vert _{H^{k_{n}}\left( Q_{T}\right) }\right)
, \text{ } \\ 
\forall \left( u,p\right) ,\left( u_{1},p_{1}\right) \in \overline{B\left(
R\right) },\text{ }\forall \gamma >0,\text{ }\forall \lambda \geq \overline{
\lambda }.%
\end{array}
\right.  \label{4.7}
\end{equation}%
\emph{In particular, numbers }$C_{1}$\emph{\ and }$\lambda $\emph{\ are also
involved in the term }$( 1+C_{1}/\lambda ^{2}) J_{2,\lambda }\left(
u,p,F_{2}\right) $\emph{\ in (\ref{4.5}). Both numbers }$\overline{ \lambda }
$\emph{\ and }$C_{1}$\emph{\ depend only on listed parameters.}

\emph{3. For every }$\lambda \geq \overline{\lambda }$\emph{\ and for every }%
$\gamma >0$\emph{\ there exists unique minimizer }$(u_{\min ,\lambda ,\gamma
}$, $p_{\min ,\lambda ,\gamma })\in \overline{B(R)}$ \emph{\ of the
functional }$J_{\lambda ,\alpha }\left( u,p\right) $\emph{\ on the set }$%
\overline{B\left( R\right) }$\emph{\ and the following inequality holds:}%
\begin{equation}
\left[ J_{\lambda ,\gamma }^{\prime }\left( u_{\min ,\lambda ,\gamma
},p_{\min ,\lambda ,\gamma }\right) ,\left( u_{\min ,\lambda ,\gamma
}-u,p_{\min ,\lambda ,\gamma }-p\right) \right] \leq 0,\text{ }\forall
\left( u,p\right) \in \overline{B\left( R\right) }.  \label{4.8}
\end{equation}

Everywhere below $C_{1}=C_{1}\left( \Omega ,T,M,R\right) >0$ \ denotes
different numbers depending only on listed parameters.

\textbf{Remark 4.1. }\emph{Even though the requirement of this theorem, so
as of all other similar theorems about the convexification method, is that
the parameter }$\lambda $\emph{\ should be sufficiently large, our
computational practice of all our previous publications about the
convexification method demonstrates that reasonable values of }$\lambda \in %
\left[ 1,5\right] $\emph{\ can always be chosen to obtain accurate
computational results, see, e.g. \cite{Bak,KL,SAR,Ktransp} and references
cited therein. On the other hand, small values of }$\lambda $\emph{\ usually
do not lead to accurate numerical results, see, e.g. \cite[Figure 1]{Ktransp}
and Table 1 in section} 5.

\textbf{Proof of Theorem 4.1.} Let $\left( u,p\right) ,\left(
u_{1},p_{1}\right) \in \overline{B\left( R\right) }$ be two arbitrary pairs
of functions. Denote $h=u_{1}-u,q=p_{1}-p.$ Hence, $\left( h,q\right) \in 
\widetilde{H}$ and%
\begin{equation}
\left. J_{\lambda ,\gamma }\left( u_{1},p_{1}\right) -J_{\lambda ,\gamma
}\left( u,p\right) =J_{\lambda ,\gamma }\left( u+h,p+q\right) -J_{\lambda
,\gamma }\left( u,p\right) .\right.  \label{4.13}
\end{equation}%
By the triangle inequality, (\ref{4.4}) and (\ref{4.12})%
\begin{equation}
\left( h,q\right) \in \overline{B_{0}\left( 2R\right) }=\left\{ \left(
h,q\right) \in \widetilde{H}:\text{ }\left\Vert h\right\Vert
_{H^{k_{n}}\left( Q_{T}\right) },\left\Vert q\right\Vert _{H^{k_{n}}\left(
Q_{T}\right) }\leq 2R\right\} .  \label{4.14}
\end{equation}%
We work with each functional $J_{1,\lambda },J_{2,\lambda }$ separately.
Note that by (\ref{2.4}) and (\ref{3.1}) 
\begin{equation}
\int\limits_{Q_{T}}\left\vert \int\limits_{\Omega }K\left( x,y\right)
g\left( y,t\right) dy\right\vert ^{2}\varphi _{\lambda }\left( t\right)
dxdt\leq C_{1}\int\limits_{Q_{T}}g^{2}\varphi _{\lambda }\left( t\right)
dxdt,\text{ }\forall g\in L_{2}\left( Q_{T}\right) .  \label{4.140}
\end{equation}

\textbf{Step 1. Analyze} $J_{1,\lambda }.$ Represent the operator $%
L_{1}\left( u+h,p+q\right) $ as the sum of its linear $L_{1,\text{lin}}$ and
nonlinear $L_{1,\text{nonlin}}$ parts with respect to $\left( h,q\right) ,$ 
\begin{equation}
\left. 
\begin{array}{c}
L_{1}\left( u+h,p+q\right) =L_{1}\left( u,p\right) +L_{1,\text{lin}}\left(
h,q\right) +L_{1,\text{nonlin}}\left( h,q\right) , \\ 
L_{1,\text{lin}}\left( h,q\right) =h_{t}+\beta \Delta h+s\nabla h\nabla
u+\dint\limits_{\Omega }K\left( x,y\right) q\left( y,t\right) dy+fq, \\ 
\text{ }L_{1,\text{nonlin}}\left( u,p\right) =s\left( x,t\right) \left(
\nabla h\right) ^{2}/2.%
\end{array}%
\right.  \label{4.15}
\end{equation}%
Using (\ref{4.5}), (\ref{4.13}) and (\ref{4.15}) we obtain 
\begin{equation}
\left. 
\begin{array}{c}
\left( L_{1}\left( u+h,p+q\right) \right) ^{2}-\left( L_{1}\left( u,p\right)
\right) ^{2}=2L_{1,\text{lin}}\left( h,q\right) L_{1}\left( u,p\right)
+\left( L_{1,\text{lin}}\left( h,q\right) \right) ^{2} \\ 
+\left[ L_{1,\text{lin}}\left( h,q\right) +L_{1}\left( u,p\right) \right]
s\left( \nabla h\right) ^{2}+s^{2}\left( \nabla h\right) ^{4}/4.%
\end{array}%
\right.  \label{4.16}
\end{equation}%
Hence, by (\ref{4.3})-(\ref{4.12}), (\ref{4.140}), (\ref{4.15}) and (\ref%
{4.16}) 
\begin{equation}
\left. 
\begin{array}{c}
J_{1,\lambda }\left( u+h,p+q\right) -J_{1,\lambda }\left( u,p\right)
=2\int\limits_{Q_{T}}L_{1,\text{lin}}\left( h,q\right) L_{1}\left(
u,p\right) \varphi _{\lambda }dxdt+ \\ 
+\widehat{J}_{1,\lambda }\left( u,p,h,q\right) , \\ 
\hspace{-0.5cm}\widehat{J}_{1,\lambda }\left( u,p,h,q\right)
=\int\limits_{Q_{T}}\left( L_{1,\text{lin}}\left( h,q\right) \right)
^{2}\varphi _{\lambda }dxdt+\int\limits_{Q_{T}}L_{1,\text{lin}}\left(
h,q\right) s\left( \nabla h\right) ^{2}\varphi _{\lambda }dxdt \\ 
+\int\limits_{Q_{T}}L_{1}\left( u,p\right) s\left( \nabla h\right)
^{2}\varphi _{\lambda }dxdt+\int\limits_{Q_{T}}\left( s^{2}\left( \nabla
h\right) ^{4}/4\right) \varphi _{\lambda }dxdt, \\ 
\lim_{\left\Vert \left( h,q\right) \right\Vert _{\widetilde{H}}\rightarrow 0}%
\left[ \left\vert \widehat{J}_{1,\lambda }\left( u,p,h,q\right) \right\vert
/\left\Vert \left( h,q\right) \right\Vert _{\widetilde{H}}\right] =0.%
\end{array}%
\right.  \label{4.17}
\end{equation}%
Denote%
\begin{equation*}
A_{1,\lambda ,\left( u,p\right) }\left( h,q\right) =2\int\limits_{Q_{T}}L_{1,%
\text{lin}}\left( h,q\right) L_{1}\left( u,p,F_{1}\right) \varphi _{\lambda
}dxdt.
\end{equation*}%
Then $A_{1,\lambda ,\left( u,p\right) }:\widetilde{H}\rightarrow \mathbb{R}$
is a bounded linear functional acting on the pair $\left( h,q\right) .$
Hence, by Riesz theorem there exists unique element $Y\in \widetilde{H}$
such that $A_{1,\lambda ,(u,p)}(h,q)$ $=[Y,\left( h,q\right) ].$ It follows
from (\ref{4.17}) that $Y=J_{1,\lambda }^{\prime }\left( u,p\right) \in 
\widetilde{H}$ is the Fr\'{e}chet derivative of the functional $J_{1,\lambda
}$ at the point $\left( u,p\right) .$ Hence, the first line of (\ref{4.17})
can be rewritten as%
\begin{equation}
J_{1,\lambda }\left( u+h,p+q\right) -J_{1,\lambda }\left( u,p\right) -\left[
J_{1,\lambda }^{\prime }\left( u,p\right) ,\left( h,q\right) \right] =%
\widehat{J}_{1,\lambda }\left( u,p,h,q\right) .  \label{4.18}
\end{equation}

We now work with the right hand side of (\ref{4.18}). Applying
Cauchy-Schwarz inequality, we obtain%
\begin{equation}
\left. 
\begin{array}{c}
\int\limits_{Q_{T}}L_{1,\text{lin}}\left( h,q\right) s\left( \nabla h\right)
^{2}\varphi _{\lambda }dxdt\geq \\ 
\geq -\left( 1/2\right) \int\limits_{Q_{T}}\left( L_{1,\text{lin}}\left(
h,q\right) \right) ^{2}\varphi _{\lambda }dxdt-\left( 1/2\right)
\int\limits_{Q_{T}}s^{2}\left( \nabla h\right) ^{4}\varphi _{\lambda }dxdt.%
\end{array}%
\right.  \label{4.180}
\end{equation}%
By (\ref{4.3}) and (\ref{4.14}) 
\begin{equation}
\left( \nabla h\left( x,t\right) \right) ^{4}\leq C_{1}\left( \nabla h\left(
x,t\right) \right) ^{2}\text{ in }Q_{T}.  \label{4.181}
\end{equation}%
Hence, (\ref{4.140}), (\ref{4.15}), (\ref{4.17}), (\ref{4.180}) and (\ref%
{4.181}) lead to%
\begin{equation}
\left. 
\begin{array}{c}
\widehat{J}_{1,\lambda }\left( u,p,h,q\right) \geq \left( 1/2\right)
\int\limits_{Q_{T}}\left( L_{1,\text{lin}}\left( h,q\right) \right)
^{2}\varphi _{\lambda }dxdt-C_{1}\int\limits_{Q_{T}}\left( \nabla h\right)
^{2}\varphi _{\lambda }dxdt\geq \\ 
\geq \left( 1/4\right) \int\limits_{Q_{T}}\left( h_{t}+\beta \Delta h\right)
^{2}\varphi _{\lambda }dxdt-C_{1}\int\limits_{Q_{T}}\left( \nabla h\right)
^{2}\varphi _{\lambda }dxdt-C_{1}\int\limits_{Q_{T}}q^{2}\varphi _{\lambda
}dxdt.%
\end{array}%
\right.  \label{4.19}
\end{equation}%
We now use Theorem 3.1 by applying Carleman estimate (\ref{3.2}) to the
first term in the second line of (\ref{4.19}). Choose a sufficiently large 
\begin{equation}
\lambda _{1}=\lambda _{1}\left( \Omega ,T,M,R\right) >2  \label{4.20}
\end{equation}%
and let $\lambda \geq \lambda _{1}.$ Then the term with $\left( \nabla
h\right) ^{2}$ in the second line of (\ref{4.19}) is absorbed by the term
with $\lambda $ $\left( \nabla h\right) ^{2}.$ Also, since $\left(
h,q\right) \in \overline{B_{0}\left( 2R\right) },$ then (\ref{4.12}) and (%
\ref{4.14}) imply that $h\left( x,T\right) =0.$ Thus, using (\ref{4.18}), we
obtain%
\begin{equation}
\left. 
\begin{array}{c}
J_{1,\lambda }\left( u+h,p+q\right) -J_{1,\lambda }\left( u,p\right) -\left[
J_{1,\lambda }^{\prime }\left( u,p\right) ,\left( h,q\right) \right] \\ 
+C_{1}\int\limits_{Q_{T}}q^{2}\varphi _{\lambda }dxdt\geq
C_{1}\int\limits_{Q_{T}}\left( h_{t}^{2}+\left( \Delta h\right) ^{2}\right)
\varphi _{\lambda }dxdt+ \\ 
+C_{1}\int\limits_{Q_{T}}\left( \lambda \left( \nabla h\right) ^{2}+\lambda
^{2}h^{2}\right) \varphi _{\lambda }dxdt,\text{ }\forall \lambda \geq
\lambda _{1}.%
\end{array}%
\right.  \label{4.21}
\end{equation}

\textbf{Step 2. Analyze} $J_{2,\lambda }.$ Acting completely similarly with
Step 1, we obtain that there exists Fr\'{e}chet derivative $J_{2,\lambda
}^{\prime }\left( u,p\right) \in \widetilde{H}$ of $J_{2,\lambda }$ at the
point $\left( u,p\right) \in \overline{B\left( R\right) }$ and 
\begin{equation}
\left[ J_{2,\lambda }^{\prime }\left( u,p\right) ,\left( h,q\right) \right]
=2\int\limits_{Q_{T}}L_{2,\text{lin}}\left( h,q\right) L_{2}\left(
u,p\right) \varphi _{\lambda }dxdt,  \label{4.25}
\end{equation}%
where $L_{2,\text{lin}}\left( h,q\right) $ is the linear part, with respect
to $\left( h,q\right) ,$ of $L_{2}\left( u+h,p+q\right) -L_{2}\left(
u,p\right) .$ Next, similarly with (\ref{4.18}), we obtain%
\begin{equation}
J_{2,\lambda }\left( u+h,p+q\right) -J_{2,\lambda }\left( u,p\right) -\left[
J_{2,\lambda }^{\prime }\left( u,p\right) ,\left( h,q\right) \right] = 
\widehat{J}_{2,\lambda }\left( u,p,h,q\right) .  \label{4.26}
\end{equation}

To estimate the right hand side of (\ref{4.26}) from the below, we use (\ref%
{4.25}) and (\ref{4.26}), keeping in mind that the formula for $\widehat{J}%
_{2,\lambda }\left( u,p,h,q\right) $ is similar with the formula for $%
\widehat{J}_{1,\lambda }\left( u,p,h,q\right) $ in (\ref{4.17}). Thus, 
\begin{equation}
\left. 
\begin{array}{c}
\widehat{J}_{2,\lambda }\left( u,p,h,q\right) \geq \left( 1/2\right)
\int\limits_{Q_{T}}\left( q_{t}-\beta \Delta q\right) ^{2}\varphi _{\lambda
}dxdt- \\ 
-C_{1}\int\limits_{Q_{T}}\left[ \left( \nabla q\right) ^{2}+q^{2}\right]
\varphi _{\lambda }dxdt-C_{1}\int\limits_{Q_{T}}\left[ \left( \Delta
h\right) ^{2}+\left( \nabla h\right) ^{2}\right] \varphi _{\lambda }dxdt.%
\end{array}%
\right.  \label{4.27}
\end{equation}%
Let $\lambda _{0}=\lambda _{0}\left( T,\beta ,a\right) >2$ be the parameter
of Theorem 3.2. Apply Carleman estimate (\ref{3.3}) to the second term in
the first line of (\ref{4.27}). Note that it follows from (\ref{4.12}) and (%
\ref{4.14}) that $q\left( x,0\right) =q\left( x,T\right) =0.$ Choose a
sufficiently large parameter%
\begin{equation}
\lambda _{2}=\lambda _{2}\left( \Omega ,T,M,R\right) \geq \max \left(
\lambda _{1},\lambda _{0}\right) ,  \label{4.28}
\end{equation}%
where $\lambda _{1}$ was defined in (\ref{4.20}). Let $\lambda \geq \lambda
_{2}.$ Then terms in the first integral in the second line of (\ref{4.27})
will be absorbed. Hence, using (\ref{4.26}) and (\ref{4.27}), we obtain for
all $\lambda \geq \lambda _{2}$ 
\begin{equation}
\left. 
\begin{array}{c}
J_{2,\lambda }\left( u+h,p+q\right) -J_{2,\lambda }\left( u,p\right) -\left[
J_{2,\lambda }^{\prime }\left( u,p\right) ,\left( h,q\right) \right] \geq \\ 
\hspace{-1cm}\geq C_{1}\int\limits_{Q_{T}}\left[ \sqrt{\lambda }\left(
\nabla q\right) ^{2}+\lambda ^{2}q^{2}\right] \varphi _{\lambda
}dxdt-C_{1}\int\limits_{Q_{T}}\left[ \left( \Delta h\right) ^{2}+\left(
\nabla h\right) ^{2}\right] \varphi _{\lambda }dxdt.%
\end{array}%
\right.  \label{4.29}
\end{equation}

\textbf{Step 3. Analyze} $J_{1,\lambda }+\left( 1/2+C_{1}/\lambda
^{2}\right) J_{2,\lambda }.$ It follows from (\ref{4.29}) that 
\begin{equation}
\left. 
\begin{array}{c}
\int\limits_{Q_{T}}q^{2}\varphi _{\lambda }dxdt\leq \left( C_{1}/\lambda
^{2}\right) \int\limits_{Q_{T}}\left[ \left( \Delta h\right) ^{2}+\left(
\nabla h\right) ^{2}\right] \varphi _{\lambda }dxdt+ \\ 
+\left( C_{1}/\lambda ^{2}\right) \left[ J_{2,\lambda }\left( u+h,p+q\right)
-J_{2,\lambda }\left( u,p\right) -J_{2,\lambda }^{\prime }\left( u,p\right)
\left( h,q\right) \right] ,\forall \lambda \geq \lambda _{2}.%
\end{array}%
\right.  \label{4.30}
\end{equation}%
Substituting (\ref{4.30}) in (\ref{4.21}), we obtain 
\begin{equation}
\left. 
\begin{array}{c}
J_{1,\lambda }\left( u+h,p+q\right) -J_{1,\lambda }\left( u,p\right) -\left[
J_{1,\lambda }^{\prime }\left( u,p\right) ,\left( h,q\right) \right] + \\ 
+\left( C_{1}/\lambda ^{2}\right) \left[ J_{2,\lambda }\left( u+h,p+q\right)
-J_{2,\lambda }\left( u,p\right) -\left[ J_{2,\lambda }^{\prime }\left(
u,p\right) ,\left( h,q\right) \right] \right] \geq \\ 
\hspace{-1cm}\geq C_{1}\int\limits_{Q_{T}}\left( h_{t}^{2}+\left( \Delta
h\right) ^{2}\right) \varphi _{\lambda }dxdt+C_{1}\int\limits_{Q_{T}}\left(
\lambda \left( \nabla h\right) ^{2}+\lambda ^{2}h^{2}\right) \varphi
_{\lambda }dxdt,\text{ }\forall \lambda \geq \lambda _{2}.%
\end{array}%
\right.  \label{4.31}
\end{equation}%
On the other hand, we obtain from (\ref{4.29})%
\begin{equation}
\left. 
\begin{array}{c}
C_{1}\int\limits_{Q_{T}}\left[ \left( \Delta h\right) ^{2}+\left( \nabla
h\right) ^{2}\right] \varphi _{\lambda }dxdt+ \\ 
+J_{2,\lambda }\left( u+h,p+q\right) -J_{2,\lambda }\left( u,p\right)
-J_{2,\lambda }^{\prime }\left( u,p\right) \left( h,q\right) \geq \\ 
\geq C_{1}\int\limits_{Q_{T}}\left[ \sqrt{\lambda }\left( \nabla q\right)
^{2}+\lambda ^{2}q^{2}\right] \varphi _{\lambda }dxdt,\text{ }\forall
\lambda \geq \lambda _{2}.%
\end{array}%
\right.  \label{4.310}
\end{equation}%
By (\ref{4.31}) the first line of (\ref{4.310}) can be estimated as%
\begin{equation*}
\left. 
\begin{array}{c}
C_{1}\int\limits_{Q_{T}}\left[ \left( \Delta h\right) ^{2}+\left( \nabla
h\right) ^{2}\right] \varphi _{\lambda }dxdt\leq \\ 
\leq J_{1,\lambda }\left( u+h,p+q,F_{1}\right) -J_{1,\lambda }\left(
u,p,F_{1}\right) -\left[ J_{1,\lambda }^{\prime }\left( u,p\right) ,\left(
h,q\right) \right] + \\ 
+\left( C_{1}/\lambda ^{2}\right) \left[ J_{2,\lambda }\left(
u+h,p+q,F_{2}\right) -J_{2,\lambda }\left( u,p,F_{2}\right) -\left[
J_{2,\lambda }^{\prime }\left( u,p\right) ,\left( h,q\right) \right] \right]
.%
\end{array}%
\right.
\end{equation*}%
Substituting this in (\ref{4.310}), we obtain%
\begin{equation*}
\left. 
\begin{array}{c}
J_{1,\lambda }\left( u+h,p+q\right) -J_{1,\lambda }\left( u,p\right) -\left[
J_{1,\lambda }^{\prime }\left( u,p\right) ,\left( h,q\right) \right] + \\ 
+\left( 1+C_{1}/\lambda ^{2}\right) \left[ J_{2,\lambda }\left(
u+h,p+q\right) -J_{2,\lambda }\left( u,p\right) -\left[ J_{2,\lambda
}^{\prime }\left( u,p\right) ,\left( h,q\right) \right] \right] \geq \\ 
\geq C_{1}\int\limits_{Q_{T}}\left[ \sqrt{\lambda }\left( \nabla q\right)
^{2}+\lambda ^{2}q^{2}\right] \varphi _{\lambda }dxdt,\text{ }\forall
\lambda \geq \lambda _{2}.%
\end{array}%
\right.
\end{equation*}%
Summing up this inequality with (\ref{4.31}), we obtain 
\begin{equation}
\left. 
\begin{array}{c}
J_{1,\lambda }\left( u+h,p+q\right) -J_{1,\lambda }\left( u,p\right)
-J_{1,\lambda }^{\prime }\left( u,p\right) \left( h,q\right) + \\ 
+\left( 1/2+C_{1}/\lambda ^{2}\right) \left[ J_{2,\lambda }\left(
u+h,p+q\right) -J_{2,\lambda }\left( u,p\right) -\left[ J_{2,\lambda
}^{\prime }\left( u,p\right) ,\left( h,q\right) \right] \right] \geq \\ 
\geq C_{1}\int\limits_{Q_{T}}\left( h_{t}^{2}+\left( \Delta h\right)
^{2}\right) \varphi _{\lambda }dxdt+ \\ 
+C_{1}\int\limits_{Q_{T}}\left( \lambda \left( \nabla h\right) ^{2}+\lambda
^{2}h^{2}+\sqrt{\lambda }\left( \nabla q\right) ^{2}+\lambda
^{2}q^{2}\right) \varphi _{\lambda }dxdt,\text{ }\forall \lambda \geq
\lambda _{2}.%
\end{array}%
\right.  \label{4.32}
\end{equation}%
Now, by (\ref{3.1}) $\varphi _{\lambda }\left( t\right) \geq e^{2a^{\lambda
}}$ for $t\in \left[ 0,T\right] .$ Hence, replacing in integrals in the
right hand side of (\ref{4.32}) $\varphi _{\lambda }\left( t\right) $ with $%
e^{2a^{\lambda }}$ and using (\ref{4.20}) and (\ref{4.28}), we obtain 
\begin{equation}
\left. 
\begin{array}{c}
J_{1,\lambda }\left( u+h,p+q\right) -J_{1,\lambda }\left( u,p\right) -\left[
J_{1,\lambda }^{\prime }\left( u,p\right) ,\left( h,q\right) \right] + \\ 
+\left( 1/2+C_{1}/\lambda ^{2}\right) \left[ J_{2,\lambda }\left(
u+h,p+q\right) -J_{2,\lambda }\left( u,p\right) -\left[ J_{2,\lambda
}^{\prime }\left( u,p\right) ,\left( h,q\right) \right] \right] \geq \\ 
\geq C_{1}e^{2a^{\lambda }}\left( \left\Vert \Delta h\right\Vert
_{L_{2}\left( Q_{T}\right) }^{2}+\left\Vert h\right\Vert _{H^{1,1}\left(
Q_{T}\right) }^{2}+\left\Vert q\right\Vert _{H^{1,0}\left( Q_{T}\right)
}^{2}\right) ,\text{ }\forall \lambda \geq \lambda _{2}.%
\end{array}%
\right.  \label{4.33}
\end{equation}

As to the functional $J_{3}\left( u,p\right) ,$ it obviously follows from (%
\ref{4.5})\emph{\ }that it has the Fr\'{e}chet derivative $J_{3}^{\prime
}\left( u,p\right) =2\left( u,p\right) \in \widetilde{H}$ and 
\begin{equation}
J_{3}\left( u+h,p+q\right) -J_{3}\left( u,p\right) -\left[ 2\left(
u,p\right) ,\left( h.q\right) \right] =\gamma \left( \left\Vert h\right\Vert
_{H^{k_{n}}\left( Q_{T}\right) }^{2}+\left\Vert q\right\Vert
_{H^{k_{n}}\left( Q_{T}\right) }^{2}\right) .  \label{4.34}
\end{equation}

Thus, summing up (\ref{4.33}) and (\ref{4.34}), setting $\overline{\lambda }%
=\lambda _{2}$ and taking into account the last line of (\ref{4.5}), we
obtain that there exists Fr\'{e}chet derivative $J_{\lambda ,\gamma
}^{\prime }\left( u,p\right) \in \widetilde{H}$ of the functional $%
J_{\lambda ,\gamma }\left( u,p\right) ,$ and we also obtain the following
equivalent of the target estimate (\ref{4.7}): 
\begin{equation*}
\left. 
\begin{array}{c}
J_{\lambda ,\gamma }\left( u+h,p+q\right) -J_{\lambda ,\gamma }\left(
u,p\right) -\left[ J_{\lambda ,\gamma }^{\prime }\left( u,p\right) \left(
h,q\right) \right] \geq \\ 
\geq C_{1}e^{2a^{\lambda }}\left( \left\Vert \Delta h\right\Vert
_{L_{2}\left( Q_{T}\right) }^{2}+\left\Vert h\right\Vert _{H^{1,1}\left(
Q_{T}\right) }^{2}+\left\Vert q\right\Vert _{H^{1,0}\left( Q_{T}\right)
}^{2}\right) + \\ 
+\gamma \left( \left\Vert h\right\Vert _{H^{k_{n}}\left( Q_{T}\right)
}^{2}+\left\Vert q\right\Vert _{H^{k_{n}}\left( Q_{T}\right) }^{2}\right) , 
\text{ }\forall \lambda \geq \overline{\lambda },\forall \gamma >0,%
\end{array}
\right.
\end{equation*}

Estimate (\ref{4.6}) is proved quite similarly with the proof of Theorem 3.1
of \cite{Bak}. The existence and uniqueness of the minimizer $\left( u_{\min
,\lambda ,\gamma },p_{\min ,\lambda ,\gamma }\right) \in \overline{B\left(
R\right) }$\emph{\ }of the functional $J_{\lambda ,\alpha }\left( u,p\right) 
$ on the set $\overline{B\left( R\right) }$ as well as inequality (\ref{4.8}%
) easily follow from a combination of Lemma 2.1 and Theorem 2.1 of \cite{Bak}%
. $\square $

\subsection{The accuracy of the minimizer}

\label{sec:4.3}

We now use the framework of the theory of Ill-Posed and Inverse Problems,
see Remark 1.1. Thus, we assume that there exist exact, noiseless data (\ref%
{2.3}), (\ref{2.5}), which are generated by the exact solution $\left(
u^{\ast },p^{\ast }\right) \left( x,t\right) $ of the MFGS (\ref{2.1}), 
\begin{equation}
u^{\ast }\left( x,T\right) =u_{T}^{\ast }\left( x\right) ,p^{\ast }\left(
x,0\right) =p_{0}^{\ast }\left( x\right) ,p^{\ast }\left( x,T\right)
=p_{T}^{\ast }\left( x\right) ,\text{ }x\in \Omega .  \label{4.35}
\end{equation}%
Let a sufficiently small number $\delta \in \left( 0,1\right) $ be the level
of the noise in the data. We assume that, similarly with (\ref{4.4}), 
\begin{equation}
\left( u^{\ast },p^{\ast }\right) \in B^{\ast }\left( R,\delta \right)
=\left\{ 
\begin{array}{c}
\left( u,p\right) \in H^{k_{n}}\left( Q_{T}\right) \times H^{k_{n}}\left(
Q_{T}\right) : \\ 
\text{ }u,p\in H_{0}^{2}\left( Q_{T}\right) , \\ 
u\left( x,T\right) =u_{T}^{\ast }\left( x\right) ,\text{ }p\left( x,T\right)
=p_{T}^{\ast }\left( x\right) , \\ 
\text{ }p\left( x,0\right) =p_{0}^{\ast }\left( x\right) , \\ 
\left\Vert u\right\Vert _{H^{k_{n}}\left( Q_{T}\right) },\left\Vert
p\right\Vert _{H^{k_{n}}\left( Q_{T}\right) }<R, \\ 
\left\Vert \Delta u\right\Vert _{L_{2}\left( Q_{T}\right) }+\left\Vert
u\right\Vert _{H^{1,1}\left( Q_{T}\right) }+ \\ 
+\left\Vert p\right\Vert _{H^{1,0}\left( Q_{T}\right) }\leq R-C_{1}\delta ,
\\ 
\left\Vert u_{T}^{\ast }\right\Vert _{H^{k_{n}}\left( \Omega \right)
}+\left\Vert p_{T}^{\ast }\right\Vert _{H^{k_{n}}\left( \Omega \right) }+ \\ 
+\left\Vert p_{0}^{\ast }\right\Vert _{H^{k_{n}}\left( \Omega \right)
},<R-C_{1}\delta .%
\end{array}%
\right\}  \label{4.36}
\end{equation}%
Since $\delta $ is sufficiently small, then it is reasonable to assume that $%
R-C_{1}\delta >0$ in (\ref{4.36}). Consider noisy data (\ref{2.3}), (\ref%
{2.5}),%
\begin{equation}
u\left( x,T\right) =u_{T}\left( x\right) ,\text{ }p\left( x,0\right)
=p_{0}\left( x\right) ,\text{ }p\left( x,T\right) =p_{T}\left( x\right) ,%
\text{ }x\in \Omega .  \label{4.37}
\end{equation}%
If the noise is random and non-smooth, as it always the case is in practice,
one can always smooth it out by one of the well known methods, see, e.g. 
\cite[Test 3]{Kpar}. Thus, we assume that 
\begin{equation}
\left. 
\begin{array}{c}
u_{T}^{\ast },\text{ }p_{0}^{\ast },\text{ }p_{T}^{\ast },\text{ }u_{T},%
\text{ }p_{0},\text{ }p_{T}\in H^{k_{n}}\left( \Omega \right) , \\ 
\left\Vert u_{T}-u_{T}^{\ast }\right\Vert _{H^{k_{n}}\left( \Omega \right) },%
\text{ }\left\Vert p_{0}-p_{0}^{\ast }\right\Vert _{H^{k_{n}}\left( \Omega
\right) },\left\Vert p_{T}-p_{T}^{\ast }\right\Vert _{H^{k_{n}}\left( \Omega
\right) }<C_{1}\delta , \\ 
\partial _{n}u_{T}\mid _{\partial \Omega }=\partial _{n}p_{0}\mid _{\partial
\Omega }=\partial _{n}p_{T}\mid _{\partial \Omega }=0.%
\end{array}%
\right.  \label{4.38}
\end{equation}

Denote 
\begin{equation}
\left. 
\begin{array}{c}
g_{1}^{\ast }\left( x,t\right) =\left( t/T\right) u_{T}^{\ast }\left(
x\right) ,\text{ }g_{2}^{\ast }\left( x,t\right) =\left( t/T\right)
p_{T}^{\ast }\left( x\right) +\left( 1-t/T\right) p_{0}^{\ast }\left(
x\right) , \\ 
\text{ }g_{1}\left( x,t\right) =\left( t/T\right) u_{T}\left( x\right) ,%
\text{ }g_{2}\left( x,t\right) =\left( t/T\right) p_{T}\left( x\right)
+\left( 1-t/T\right) p_{0}\left( x\right) , \\ 
\widetilde{u}^{\ast }=u^{\ast }-g_{1}^{\ast },\text{ }\widetilde{p}^{\ast
}=p^{\ast }-g_{2}^{\ast }, \\ 
\widetilde{u}=u-g_{1},\text{ }\widetilde{p}=p-g_{2},\text{ }\forall \left(
u,p\right) \in B\left( R\right) .%
\end{array}%
\right.  \label{4.40}
\end{equation}%
Using the second line of (\ref{4.38}), we obtain 
\begin{equation}
\left\Vert g_{1}-g_{1}^{\ast }\right\Vert _{H^{k_{n}}\left( Q_{T}\right)
},\left\Vert g_{2}-g_{2}^{\ast }\right\Vert _{H^{k_{n}}\left( Q_{T}\right)
}<C_{1}\delta .  \label{4.41}
\end{equation}%
Then by (\ref{4.4}), (\ref{4.12}), (\ref{4.14}), (\ref{4.36}), (\ref{4.40})
and (\ref{4.41})%
\begin{equation}
\left( \widetilde{u}^{\ast },\widetilde{p}^{\ast }\right) ,\text{ }\left( 
\widetilde{u},\widetilde{p}\right) \in B_{0}\left( 2R\right) .  \label{4.43}
\end{equation}%
Also, by (\ref{2.1}) and (\ref{4.35}) 
\begin{equation}
L_{1}\left( \widetilde{u}^{\ast }+g_{1}^{\ast },\text{ }\widetilde{p}^{\ast
}+g_{2}^{\ast }\right) =0,\text{ }L_{2}\left( \widetilde{u}^{\ast
}+g_{1}^{\ast },\text{ }\widetilde{p}^{\ast }+g_{2}^{\ast }\right) =0.
\label{4.44}
\end{equation}%
Based on (\ref{4.40}) and (\ref{4.43}), consider a new functional $%
I_{\lambda ,\gamma },$%
\begin{equation}
I_{\overline{\lambda }_{2},\gamma }:B_{0}\left( 2R\right) \rightarrow 
\mathbb{R},\text{ }I_{\overline{\lambda }_{2},\gamma }\left( v,w\right) =J_{%
\overline{\lambda }_{2},\gamma }\left( v+g_{1},w+g_{2}\right) .  \label{4.45}
\end{equation}

\textbf{Theorem \ 4.2}. \emph{Let conditions of Theorem 4.1 as well as
conditions (\ref{4.35})-(\ref{4.38}) hold and }$\delta \in \left( 0,1\right) 
$ \emph{is so small that in (\ref{4.36}) }$R-C_{1}\delta >0$\emph{. Let }$%
\overline{\lambda }=\overline{\lambda }\left( \Omega ,T,M,R\right) >1$\emph{%
\ \ be the number chosen in Theorem 4.1. Consider the number }$\overline{%
\lambda }_{2}=\overline{\lambda }\left( \Omega ,T,M,2R\right) .$\emph{\
Without any loss of generality, we assume that }$\overline{\lambda }_{2}\geq 
\overline{\lambda }.$ \emph{As it is often done in the regularization theory 
\cite{T}, choose the regularization parameter }$\gamma $\emph{\ depending on
the noise level }$\delta $\emph{\ as }$\gamma \left( \delta \right) =\delta
^{2}.$ \emph{Then:}

1.\emph{\ There exists unique minimizer }$\left( v_{\min ,\overline{\lambda }%
_{2},\gamma },w_{\min ,\overline{\lambda }_{2},\gamma }\right) \in \overline{%
B_{0}\left( 2R\right) }$\emph{\ of the functional} $I_{\overline{\lambda }%
_{2},\gamma }\left( v,w\right) $\emph{\ on the set }$\overline{B_{0}\left(
2R\right) }.$ \emph{The following accuracy estimates hold:}%
\begin{equation}
\left. 
\begin{array}{c}
\left\Vert \Delta \widetilde{u}^{\ast }-\Delta v_{\min ,\overline{\lambda }%
_{2},\gamma \left( \delta \right) }\right\Vert _{L_{2}\left( Q_{T}\right)
}+\left\Vert \widetilde{u}^{\ast }-v_{\min ,\overline{\lambda }_{2},\gamma
\left( \delta \right) }\right\Vert _{H^{1,1}\left( Q_{T}\right) }+ \\ 
+\left\Vert \widetilde{p}^{\ast }-w_{\min ,\overline{\lambda }_{2},\gamma
\left( \delta \right) }\right\Vert _{H^{1,0}\left( Q_{T}\right) }\leq
C_{1}\delta .%
\end{array}%
\right.  \label{4.51}
\end{equation}

2. \emph{Define functions }$\overline{u}_{\min ,\overline{\lambda }
_{2},\gamma \left( \delta \right) }$ \emph{and} $\overline{p}_{\min , 
\overline{\lambda }_{2}\gamma \left( \delta \right) },$ \emph{\ }%
\begin{equation}
\overline{u}_{\min ,\overline{\lambda }_{2},\gamma \left( \delta \right)
}=v_{\min ,\overline{\lambda }_{2},\gamma \left( \delta \right) }+g_{1}, 
\text{ }\overline{p}_{\min ,\overline{\lambda }_{2}\gamma \left( \delta
\right) }=w_{\min ,\overline{\lambda }_{2},\gamma \left( \delta \right)
}+g_{2}.  \label{4.450}
\end{equation}%
\emph{Then } \emph{\ }%
\begin{equation}
\left\Vert \Delta \overline{u}_{\min ,\overline{\lambda }_{2},\gamma \left(
\delta \right) }\right\Vert _{L_{2}\left( Q_{T}\right) },\left\Vert 
\overline{u}_{\min ,\overline{\lambda }_{2},\gamma \left( \delta \right)
}\right\Vert _{H^{1,1}\left( Q_{T}\right) },\left\Vert \overline{p}_{\min , 
\overline{\lambda }_{2},\gamma \left( \delta \right) }\right\Vert
_{H^{1,1}\left( Q_{T}\right) }\leq R  \label{4.451}
\end{equation}%
\emph{\ and the following accuracy estimates hold:}%
\begin{equation}
\left. 
\begin{array}{c}
\left\Vert \Delta u^{\ast }-\Delta \overline{u}_{\min ,\overline{\lambda }
_{2},\gamma \left( \delta \right) }\right\Vert _{L_{2}\left( Q_{T}\right)
}+\left\Vert u^{\ast }-\overline{u}_{\min ,\overline{\lambda }_{2},\gamma
\left( \delta \right) }\right\Vert _{H^{1,1}\left( Q_{T}\right) }+ \\ 
+\left\Vert p^{\ast }-\overline{p}_{\min ,\overline{\lambda }_{2}\gamma
\left( \delta \right) }\right\Vert _{H^{1,0}\left( Q_{T}\right) }\leq
C_{1}\delta .%
\end{array}
\right.  \label{4.39}
\end{equation}

3. \emph{Next, based on (\ref{4.4}), (\ref{4.36}), (\ref{4.450})-(\ref{4.39}%
), assume that }%
\begin{equation}
\left( \overline{u}_{\min ,\overline{\lambda }_{2},\gamma \left( \delta
\right) },\overline{p}_{\min ,\overline{\lambda }_{2},\gamma \left( \delta
\right) }\right) \in \overline{B\left( R\right) }.  \label{4.390}
\end{equation}%
\emph{\ Then this vector function is the unique minimizer of the functional }%
$J_{\overline{\lambda }_{2},\gamma }\left( u,p\right) $\emph{\ on the set }$%
\overline{B\left( R\right) },$\emph{\ which is found in Theorem 4.1, i.e. }%
\begin{equation}
\left( \overline{u}_{\min ,\overline{\lambda }_{2},\gamma \left( \delta
\right) },\overline{p}_{\min ,\overline{\lambda }_{2},\gamma \left( \delta
\right) }\right) =\left( u_{\min ,\overline{\lambda }_{2},\gamma \left(
\delta \right) },p_{\min ,\overline{\lambda }_{2},\gamma \left( \delta
\right) }\right) .  \label{4.391}
\end{equation}%
\emph{Therefore, accuracy estimates (\ref{4.39}) remain valid if the pair }$%
\big(\overline{u}_{\min ,\overline{\lambda }_{2},\gamma \left( \delta
\right) },\overline{p}_{\min ,\overline{\lambda }_{2},\gamma \left( \delta
\right) }\big)$\emph{\ is replaced with the pair }$\left( u_{\min ,\overline{%
\lambda }_{2},\gamma \left( \delta \right) },p_{\min ,\overline{\lambda }%
_{2},\gamma \left( \delta \right) }\right) .$

\textbf{Proof}. We note first that inequalities (\ref{4.451}) are
insufficient to ensure (\ref{4.390}). This is because the norm $\left\Vert
.\right\Vert _{H^{k_{n}}\left( Q_{T}\right) }$ is stronger than the norms
presented in (\ref{4.451}), see (\ref{4.01}). The full analog of Theorem 4.1
is obviously valid for the functional $I_{\lambda ,\gamma }\left( v,w\right) 
$ with the only change of $\overline{\lambda }$ with $\overline{\lambda }%
_{2}\geq \overline{\lambda }$ and $B\left( R\right) $ with $B_{0}\left(
2R\right) .$ Hence, the existence of the unique minimizer $\left( v_{\min ,%
\overline{\lambda }_{2},\gamma },w_{\min ,\overline{\lambda }_{2},\gamma
}\right) \in \overline{B_{0}\left( 2R\right) }$ of the functional $I_{%
\overline{\lambda }_{2},\gamma }\left( v,w\right) $ on the set $\overline{%
B_{0}\left( 2R\right) }$ follows from Theorem 4.1. By (\ref{4.40}) and (\ref%
{4.45}) 
\begin{equation}
I_{\overline{\lambda }_{2},\gamma }\left( \widetilde{u}^{\ast },\widetilde{p}%
^{\ast }\right) =J_{\overline{\lambda }_{2},\gamma }\left( \widetilde{u}%
^{\ast }+g_{1},\widetilde{p}^{\ast }+g_{2}\right) .  \label{4.46}
\end{equation}%
Let $I_{\overline{\lambda }_{2},\gamma }^{\prime }\left( v,w\right) $ be the
Fr\'{e}chet derivative of the functional $I_{\overline{\lambda }_{2},\gamma
} $ at the point $\left( v,w\right) .$ By (\ref{4.7}), (\ref{4.40}), (\ref%
{4.43}) and (\ref{4.46})%
\begin{equation}
\left. 
\begin{array}{c}
I_{\overline{\lambda }_{2},\gamma }\left( \widetilde{u}^{\ast },\widetilde{p}%
^{\ast }\right) -I_{\overline{\lambda }_{2},\gamma }\left( v_{\min ,%
\overline{\lambda }_{2},\gamma },w_{\min ,\overline{\lambda }_{2},\gamma
}\right) - \\ 
-\left[ I_{\overline{\lambda }_{2},\gamma }^{\prime }\left( v_{\min ,%
\overline{\lambda }_{2},\gamma },w_{\min ,\overline{\lambda }_{2},\gamma
}\right) ,\widetilde{u}^{\ast }-v_{\min ,\overline{\lambda }_{2},\gamma },%
\widetilde{p}^{\ast }-w_{\min ,\overline{\lambda }_{2},\gamma }\right] \geq
\\ 
\geq C_{1}\left( \left\Vert \Delta \widetilde{u}^{\ast }-\Delta v_{\min ,%
\overline{\lambda }_{2},\gamma }\right\Vert _{L_{2}\left( Q_{T}\right)
}^{2}+\left\Vert \widetilde{u}^{\ast }-v_{\min ,\overline{\lambda }%
_{2},\gamma }\right\Vert _{H^{1,1}\left( Q_{T}\right) }^{2}\right) + \\ 
+C_{1}\left\Vert \widetilde{p}^{\ast }-w_{\min ,\overline{\lambda }%
_{2},\gamma }\right\Vert _{H^{1,0}\left( Q_{T}\right) }^{2}.%
\end{array}%
\right.  \label{4.47}
\end{equation}%
Consider first the left hand side of inequality (\ref{4.47}). By (\ref{4.8}) 
\begin{equation*}
-\left[ I_{\lambda ,\gamma }^{\prime }\left( v_{\min ,\overline{\lambda }%
_{2},\gamma },w_{\min ,\overline{\lambda }_{2},\gamma }\right) ,\widetilde{u}%
^{\ast }-v_{\min ,\overline{\lambda }_{2},\gamma },\widetilde{p}^{\ast
}-w_{\min ,\overline{\lambda }_{2},\gamma }\right] \leq 0.
\end{equation*}%
Also, $-I_{\lambda ,\gamma }\left( v_{\min ,\overline{\lambda }_{2},\gamma
},w_{\min ,\overline{\lambda }_{2},\gamma }\right) \leq 0.$ Hence, the left
hand side of (\ref{4.47}) can be estimated as: 
\begin{equation}
\left. 
\begin{array}{c}
I_{\overline{\lambda }_{2},\gamma }\left( \widetilde{u}^{\ast },\widetilde{p}%
^{\ast }\right) -I_{\overline{\lambda }_{2},\gamma }\left( v_{\min ,%
\overline{\lambda }_{2},\gamma },w_{\min ,\overline{\lambda }_{2},\gamma
}\right) - \\ 
-\left[ I_{\overline{\lambda }_{2},\gamma }^{\prime }\left( v_{\min ,%
\overline{\lambda }_{2},\gamma },w_{\min ,\overline{\lambda }_{2},\gamma
}\right) ,\left( \widetilde{u}^{\ast }-v_{\min ,\overline{\lambda }%
_{2},\gamma },\widetilde{p}^{\ast }-w_{\min ,\overline{\lambda }_{2},\gamma
}\right) \right] \leq \\ 
\leq I_{\overline{\lambda }_{2},\gamma }\left( \widetilde{u}^{\ast },%
\widetilde{p}^{\ast }\right) .%
\end{array}%
\right.  \label{4.48}
\end{equation}

We now estimate from the above the right hand side of (\ref{4.48}). By (\ref%
{4.5}) and (\ref{4.46})%
\begin{equation}
\left. 
\begin{array}{c}
I_{\overline{\lambda }_{2},\gamma }\left( \widetilde{u}^{\ast },\widetilde{p}
^{\ast }\right) =J_{\overline{\lambda }_{2},\gamma }\left( \left( \widetilde{
u}^{\ast }+g_{1}^{\ast }\right) +\left( g_{1}-g_{1}^{\ast }\right) ,\left( 
\widetilde{p}^{\ast }+g_{2}^{\ast }\right) +\left( g_{2}-g_{2}^{\ast
}\right) \right) = \\ 
=\int\limits_{Q_{T}}\left[ L_{1}\left( \left( \widetilde{u}^{\ast
}+g_{1}^{\ast }\right) +\left( g_{1}-g_{1}^{\ast }\right) \right) \right]
^{2}\varphi _{\overline{\lambda }_{2}}dxdt+ \\ 
+\left( 1/2+C_{1}/\overline{\lambda }_{2}^{2}\right) \int\limits_{Q_{T}} %
\left[ L_{2}\left( \left( \widetilde{u}^{\ast }+g_{1}^{\ast }\right) +\left(
g_{1}-g_{1}^{\ast }\right) \right) \right] ^{2}\varphi _{\overline{\lambda }
_{2}}dxdt+ \\ 
+\gamma \left( \left\Vert \widetilde{u}^{\ast }+g_{1}\right\Vert
_{H^{k_{n}}\left( Q_{T}\right) }^{2}+\left\Vert \widetilde{p}^{\ast
}+g_{2}\right\Vert _{H^{k_{n}}\left( Q_{T}\right) }^{2}\right) .%
\end{array}
\right.  \label{4.49}
\end{equation}%
By (\ref{3.1}) $\varphi _{\lambda }\left( t\right) \leq e^{2\left(
T+a\right) ^{\lambda }}$ for $t\in \left[ 0,T\right] .$ Hence, using (\ref%
{2.1}), (\ref{2.4}), (\ref{4.140}), (\ref{4.41}) and (\ref{4.44}), we obtain%
\begin{equation*}
\int\limits_{Q_{T}}\left[ L_{i}\left( \left( \widetilde{u}^{\ast
}+g_{1}^{\ast }\right) +\left( g_{1}-g_{1}^{\ast }\right) \right) \right]
^{2}\varphi _{\lambda }dxdt\leq C_{1}e^{2\left( T+a\right) ^{\lambda
}}\delta ^{2},\text{ }i=1,2.
\end{equation*}%
Hence, (\ref{4.43}) and (\ref{4.49}) lead to 
\begin{equation}
I_{\overline{\lambda }_{2},\gamma }\left( \widetilde{u}^{\ast },\widetilde{p}
^{\ast }\right) \leq C_{1}e^{2\left( T+a\right) ^{\overline{\lambda }
_{2}}}\delta ^{2}+C_{1}\gamma .  \label{4.50}
\end{equation}%
Since $\overline{\lambda }_{2}=\overline{\lambda }\left( \left( \Omega
,T,M,2R\right) \right) ,$ then we can set the right hand side of (\ref{4.50}%
) as $C_{1}\delta ^{2}+C_{1}\gamma .$ Setting $\gamma =\gamma \left( \delta
\right) =\delta ^{2}$ and using (\ref{4.47}), (\ref{4.48}) and (\ref{4.50}),
we obtain (\ref{4.51}), which is our first target estimate.

We now use (\ref{4.40}) and (\ref{4.450}) to prove (\ref{4.451}) and (\ref%
{4.39}). We have 
\begin{equation}
\left. 
\begin{array}{c}
\widetilde{u}^{\ast }-v_{\min ,\overline{\lambda }_{2},\gamma \left( \delta
\right) }=\left( u^{\ast }-\overline{u}_{\min ,\overline{\lambda }%
_{2},\gamma \left( \delta \right) }\right) -\left( g_{1}^{\ast
}-g_{1}\right) , \\ 
\widetilde{p}^{\ast }-w_{\min ,\overline{\lambda }_{2},\gamma \left( \delta
\right) }=\left( p^{\ast }-\overline{p}_{\min ,\overline{\lambda }%
_{2},\gamma \left( \delta \right) }\right) -\left( g_{2}^{\ast
}-g_{2}\right) .%
\end{array}%
\right.   \label{4.500}
\end{equation}%
Hence, by (\ref{4.41}), (\ref{4.500}) and triangle inequality%
\begin{equation*}
\left. 
\begin{array}{c}
\left\Vert \Delta \widetilde{u}^{\ast }-\Delta v_{\min ,\overline{\lambda }%
_{2},\gamma \left( \delta \right) }\right\Vert _{L_{2}\left( Q_{T}\right)
}+\left\Vert \widetilde{u}^{\ast }-v_{\min ,\overline{\lambda }_{2},\gamma
\left( \delta \right) }\right\Vert _{H^{1,1}\left( Q_{T}\right) }+ \\ 
+\left\Vert \widetilde{p}^{\ast }-w_{\min ,\overline{\lambda }_{2},\gamma
\left( \delta \right) }\right\Vert _{H^{1,0}\left( Q_{T}\right) }\geq
\left\Vert \Delta u^{\ast }-\Delta \overline{u}_{\min ,\overline{\lambda }%
_{2},\gamma \left( \delta \right) }\right\Vert _{L_{2}\left( Q_{T}\right) }+
\\ 
+\left\Vert u^{\ast }-\overline{u}_{\min ,\overline{\lambda }_{2},\gamma
\left( \delta \right) }\right\Vert _{H^{1,1}\left( Q_{T}\right) }+\left\Vert
p^{\ast }-\overline{p}_{\min ,\overline{\lambda }_{2},\gamma \left( \delta
\right) }\right\Vert _{H^{1,0}\left( Q_{T}\right) }-C_{1}\delta .%
\end{array}%
\right. 
\end{equation*}%
Comparing this with (\ref{4.51}), we obtain (\ref{4.39}). Next, the last
line of (\ref{4.36}), (\ref{4.39}) $\ $and triangle inequality imply (\ref%
{4.451}).

Assume now that (\ref{4.390}) is valid. Consider the minimizer $\big(
u_{\min ,\overline{\lambda }_{2},\gamma \left( \delta \right) }$, $p_{\min ,%
\overline{\lambda }_{2},\gamma \left( \delta \right) } \big)$ $\in \overline{%
B\left( R\right) }$ of the functional $J_{\overline{\lambda }_{2},\gamma
\left( \delta \right) }$ on the set $\overline{B\left( R\right) },$ which
was found in Theorem 4.1. Consider the pair of functions $V_{\min ,\overline{%
\lambda }_{2},\gamma \left( \delta \right) },W_{\min ,\overline{\lambda }%
_{2},\gamma \left( \delta \right) }$ defined as: 
\begin{equation}
V_{\min ,\overline{\lambda }_{2},\gamma \left( \delta \right) }=u_{\min ,%
\overline{\lambda }_{2},\gamma \left( \delta \right) }-g_{1},\text{ }W_{\min
,\overline{\lambda }_{2},\gamma \left( \delta \right) }=p_{\min ,\overline{%
\lambda }_{2},\gamma \left( \delta \right) }-g_{2}.  \label{4.70}
\end{equation}%
Then, similarly with (\ref{4.43}), $\left( V_{\min ,\overline{\lambda }%
_{2},\gamma \left( \delta \right) },W_{\min ,\overline{\lambda }_{2},\gamma
\left( \delta \right) }\right) \in \overline{B_{0}\left( 2R\right) }.$ Since
the vector function $\left( v_{\min ,\overline{\lambda }_{2},\gamma
},w_{\min ,\overline{\lambda }_{2},\gamma }\right) $ is the unique minimizer
of the functional $I_{\overline{\lambda }_{2},\gamma \left( \delta \right) }$
on the set $\overline{B_{0}\left( 2R\right) },$ then $I_{\overline{\lambda }%
_{2},\gamma \left( \delta \right) }\left( v_{\min ,\overline{\lambda }%
_{2},\gamma },w_{\min ,\overline{\lambda }_{2},\gamma }\right) \leq I_{%
\overline{\lambda }_{2},\gamma \left( \delta \right) }\left( V_{\min ,%
\overline{\lambda }_{2},\gamma \left( \delta \right) },W_{\min ,\overline{%
\lambda }_{2},\gamma \left( \delta \right) }\right) .$ Hence, using (\ref%
{4.45}), (\ref{4.450}) and (\ref{4.70}), we obtain $J_{\overline{\lambda }%
_{2},\gamma \left( \delta \right) }\left( \overline{u}_{\min ,\overline{%
\lambda }_{2},\gamma \left( \delta \right) },\overline{p}_{\min ,\overline{%
\lambda }_{2},\gamma \left( \delta \right) }\right) \leq J_{\overline{%
\lambda }_{2},\gamma \left( \delta \right) }\left( u_{\min ,\overline{%
\lambda }_{2},\gamma \left( \delta \right) },_{\min ,\overline{\lambda }%
_{2},\gamma \left( \delta \right) }\right) .$ Hence, the vector function $(%
\overline{u}_{\min ,\overline{\lambda }_{2},\gamma (\delta )},\overline{p}%
_{\min ,\overline{\lambda }_{2},\gamma \left( \delta \right) })$ is also a
minimizer of the functional $J_{\overline{\lambda }_{2},\gamma \left( \delta
\right) }$ on the set $\overline{B\left( R\right) }.$ Since the minimizer on 
$\overline{B\left( R\right) }$ is unique by Theorem 4.1, then (\ref{4.391})
is valid. $\square $

\subsection{The global convergence of the gradient descent method}

\label{sec:4.4}

Assume now that in (\ref{4.36}) $R$ is replaced with $R/3,$ i.e. assume that%
\begin{equation}
\left( u^{\ast },p^{\ast }\right) \in B^{\ast }\left( R/3,\delta \right) 
\text{ and }R/3-C_{1}\delta >0.  \label{4.60}
\end{equation}%
The inequality in (\ref{4.60}) is reasonable since the noise level $\delta
\in \left( 0,1\right) $ is small. Suppose that conditions of Theorem 4.2
hold. Let the assumption (\ref{4.390}) be valid. Then it follows from (\ref%
{4.39}), (\ref{4.391}) and (\ref{4.60}) that it is reasonable to assume that%
\begin{equation}
\left( u_{\min ,\overline{\lambda }_{2},\gamma \left( \delta \right)
},p_{\min ,\overline{\lambda }_{2},\gamma \left( \delta \right) }\right) \in
B\left( R/3\right) .  \label{4.62}
\end{equation}

We now construct the gradient descent method of the minimization of the
functional $J_{\overline{\lambda }_{2},\gamma \left( \delta \right) }.$ Let $%
\varkappa \in \left( 0,1\right) $ be its step size. Consider an arbitrary
point $\left( u^{\left( 0\right) },p^{\left( 0\right) }\right) ,$ 
\begin{equation}
\left( u^{\left( 0\right) },p^{\left( 0\right) }\right) \in B\left(
R/3\right) .  \label{4.63}
\end{equation}%
The iterative sequence of the gradient descent method \ is:%
\begin{equation}
\left( u^{\left( n\right) },p^{\left( n\right) }\right) =\left( u^{\left(
n-1\right) },p^{\left( n-1\right) }\right) -\varkappa J_{\overline{\lambda }%
_{2},\gamma \left( \delta \right) }^{\prime }\left( u^{\left( n-1\right)
},p^{\left( n-1\right) }\right) ,\text{ }n=1,2,...  \label{4.64}
\end{equation}%
where $J_{\overline{\lambda }_{2},\gamma \left( \delta \right) }^{\prime
}\left( u^{\left( n-1\right) },p^{\left( n-1\right) }\right) \in \widetilde{H%
}$ is the Fr\'{e}chet derivative of the functional $J_{\overline{\lambda }%
_{2},\gamma \left( \delta \right) }$ at the point $\left( u^{\left(
n-1\right) },p^{\left( n-1\right) }\right) .$ Note that (\ref{4.4}), (\ref%
{4.12}), (\ref{4.63}) and (\ref{4.64}) imply that 
\begin{equation*}
u^{\left( n\right) }\left( x,T\right) =u_{T}\left( x\right) ,\text{ }%
p^{\left( n\right) }\left( x,0\right) =p_{0}\left( x\right) ,\text{ }%
p^{\left( n\right) }\left( x,T\right) =p_{T}\left( x\right) ,\text{ }%
n=1,2,...
\end{equation*}

\textbf{Theorem 4.3. }\emph{Assume that conditions of Theorem 4.2 are
satisfied as well as conditions (\ref{4.60})-(\ref{4.64}). Then there exists
a number }$\varkappa _{0}\in \left( 0,1\right) $\emph{\ such that for any }$%
\varkappa \in \left( 0,\varkappa _{0}\right) $\emph{\ there exists a number }%
$\theta =\theta \left( \varkappa \right) \in \left( 0,1\right) $\emph{\ such
that }%
\begin{equation}
\left( u^{\left( n\right) },p^{\left( n\right) }\right) \in B\left( R\right)
,\text{ }n=1,2,...  \label{4.65}
\end{equation}%
\emph{\ and the following convergence estimate is valid:}%
\begin{equation}
\left. 
\begin{array}{c}
\left\Vert \Delta u^{\ast }-\Delta u^{\left( n\right) }\right\Vert
_{L_{2}\left( Q_{T}\right) }+\left\Vert u^{\ast }-u^{\left( n\right)
}\right\Vert _{H^{1,1}\left( Q_{T}\right) }+\left\Vert p^{\ast }-p^{\left(
n\right) }\right\Vert _{H^{1,0}\left( Q_{T}\right) }\leq \\ 
\hspace{-1cm}\leq C_{1}\delta +\theta ^{n}\left( \left\Vert u_{\min ,%
\overline{\lambda }_{2},\gamma \left( \delta \right) }-u^{\left( 0\right)
}\right\Vert _{H^{k_{n}}\left( Q_{T}\right) }+\left\Vert p_{\min ,\overline{%
\lambda }_{2},\gamma \left( \delta \right) }-p^{\left( 0\right) }\right\Vert
_{H^{k_{n}}\left( Q_{T}\right) }\right) .%
\end{array}%
\right.  \label{4.66}
\end{equation}

\textbf{Proof.} The existence of numbers $\varkappa _{0}$ and $\theta \left(
\varkappa \right) $ as well as (\ref{4.65}) and estimate%
\begin{equation}
\left. 
\begin{array}{c}
\left\Vert u_{\min ,\overline{\lambda }_{2},\gamma \left( \delta \right)
}-u^{\left( n\right) }\right\Vert _{H^{k_{n}}\left( Q_{T}\right)
}+\left\Vert p_{\min ,\overline{\lambda }_{2},\gamma \left( \delta \right)
}-p^{\left( n\right) }\right\Vert _{H^{k_{n}}\left( Q_{T}\right) }\leq \\ 
\leq \theta ^{n}\left( \left\Vert u_{\min ,\overline{\lambda }_{2},\gamma
\left( \delta \right) }-u^{\left( 0\right) }\right\Vert _{H^{k_{n}}\left(
Q_{T}\right) }+\left\Vert p_{\min ,\overline{\lambda }_{2},\gamma \left(
\delta \right) }-p^{\left( 0\right) }\right\Vert _{H^{k_{n}}\left(
Q_{T}\right) }\right)%
\end{array}%
\right.  \label{4.67}
\end{equation}%
follow immediately from (\ref{4.62})-(\ref{4.64}) and \cite[Theorem 6]{SAR}.
Next, triangle inequality, (\ref{4.39})-(\ref{4.391}) and (\ref{4.67})
easily lead to (\ref{4.66}). $\ \square $

\textbf{Remark 4.1}. \emph{Since }$R>0$\emph{\ is an arbitrary number and by
(\ref{4.63}) the starting point of iterations }$\left( u^{\left( 0\right)
},p^{\left( 0\right) }\right) $\emph{\ is an arbitrary point of }$B\left(
R/3\right) ,$\emph{\ then Theorem 4.3 means the} \underline{global} \emph{\
convergence of the gradient descent method (\ref{4.64}).}

\section{Numerical Studies}

\label{sec:5}

In this section we describe our numerical studies of the Minimization
Problem formulated in subsection 4.1. First, as it is always done in
numerical studies of Ill-Posed and Inverse Problems (see Remark 1.1 in
section 1), we need to figure out how to numerically generate the data for
our problem. More precisely, we need to numerically generate such initial
and terminal conditions (\ref{2.3}), (\ref{2.5}) that the solution $\left(
u,p\right) \left( x,t\right) $ of problem (\ref{2.1})-(\ref{2.3}), (\ref{2.5}%
) would exist in the case of the absence of the noise in the data (\ref{2.3}%
), (\ref{2.5}): we will add the noise later. Then we need to
\textquotedblleft pretend" that we do not know the pair $\left( u,p\right)
\left( x,t\right) ,$ solve the Minimization Problem and then compare the
resulting computed solution $\left( u_{\text{comp}},p_{\text{comp}}\right)
\left( x,t\right) $ with $\left( u,p\right) \left( x,t\right) .$ In terms of
the theory of Ill-Posed and Inverse Problems (Remark 1.1 in section 1), $%
\left( u_{\text{comp}},p_{\text{comp}}\right) (x,t)$ is the reconstruction
of $\left( u,p\right) (x,t)$. Our special procedure for such data generation
is described in subsection 5.1.

\subsection{Numerical data generation}

\label{sec:5.1}

Consider an arbitrary function $v\in H^{k_{n}}$ $(Q_{T})$ with $\partial
_{\nu }v\mid _{S_{T}}=0,$ also, see (\ref{2.2}) and (\ref{4.3}). Let $%
m_{0}\left( x\right) \in H^{k_{n}}\left( \Omega \right) $ be another
arbitrary function of our choice, such that $\partial _{\nu }v\mid _{\Omega
}=0.$ Next, we solve the following initial boundary value problem:%
\begin{equation}
\left. 
\begin{array}{c}
m_{t}(x,t)-\beta \Delta m(x,t)+\func{div}(s(x,t)m(x,t)\nabla v(x,t))=0,\text{
}\left( x,t\right) \in Q_{T}, \\ 
m\left( x,0\right) =m_{0}\left( x\right) ,\text{ }\partial _{\nu }m\mid
_{S_{T}}=0.%
\end{array}%
\right.  \label{5.1}
\end{equation}%
Using the well known results of the classical theory of parabolic equations 
\cite{Lad}, we impose such conditions on the domain $\Omega $ and functions $%
s,v,$ which guarantee that the unique solution of problem (\ref{5.1}) is $%
m\in H^{k_{n}}\left( Q_{T}\right) $. Numerical solution of problem (\ref{5.1}%
) is elementary, and we compute it via the finite difference method. Next,
given functions $v\left( x,t\right) $ and $m\left( x,t\right) ,$ we compute
the function $F\left( x,t\right) ,$%
\begin{equation*}
F\left( x,t\right) =-\left( v_{t}+\beta \Delta v+s(\nabla v)^{2}/2\right)
(x,t)-\int\limits_{\Omega }K\left( x,y\right) m\left( y,t\right) dy-f\left(
x,t\right) m\left( x,t\right) .
\end{equation*}%
Thus, we have obtained the solution $\left( u,p\right) \equiv \left(
v,m\right) \in H^{k_{n}}\left( Q_{T}\right) \times H^{k_{n}}\left(
Q_{T}\right) $ of system (\ref{2.1}) with the zero Neumann boundary
conditions (\ref{2.2}) \ and with $F_{1}\left( x,t\right) =F\left(
x,t\right) ,$ $F_{2}\left( x,t\right) =0$. Initial and terminal conditions (%
\ref{2.3}), (\ref{2.5}) now are: $u\left( x,T\right) =v\left( x,T\right) ,$ $%
p(x,0)=m_{0}\left( x\right) ,$ $p(x,T)=m\left( x,T\right) .$ Uniqueness of
the solution $u,p\in H_{0}^{2}\left( Q_{T}\right) $ of problem (\ref{2.1})-(%
\ref{2.3}), (\ref{2.5}) was proven in \cite{MFG1}.

\subsection{Numerical testing for the Minimization Problem}

\label{sec:5.2}

We have conducted numerical studies in the 2D case. We took: 
\begin{equation}
\Omega =\left\{ x=(x_{1},x_{2}):x_{1},x_{2}\in \left( 0,1\right) \right\} ,%
\text{ }T=1.  \label{5.2}
\end{equation}

We take $s(x,t)=f\left( x,t\right) =1,K\left( x,y\right) =1,\beta =0.02$ in (%
\ref{2.1}), $a=1.01$ in (\ref{3.1}). As to the regularization parameter $%
\gamma $ in (\ref{4.5}), we found its optimal value $\gamma =0.001$ by the
trial and error procedure in Test 1 below. Now, it is hard to provide a
precise estimate of $C_{1}.$ But since we use $\lambda =2$ in our
computations (see below), then taking in (\ref{4.5}) $C_{1}=2,$ we obtain $%
\left( 1/2+C_{1}/\lambda ^{2}\right) =1,$ and this is what we use in (\ref%
{4.5}) in our computations.

For the initial conditions for the MFGS in (\ref{2.3}) and (\ref{5.1}), we
took 
\begin{equation}
p_{0}(x)=m_{0}(x)=\frac{1}{2\pi }\exp \left[ -\frac{%
(x_{1}-0.5)^{2}+(x_{2}-0.5)^{2}}{2}\right] .  \label{5.3}
\end{equation}

To solve problem (\ref{5.1}) for data generation, we have used the spatial
mesh sizes $1/80\times 1/80$ and the temporal step size $1/320$. In the
computations of the Minimization Problem, the spatial mesh sizes were $%
1/20\times 1/20.$ and the temporal step size was $1/10$.

To guarantee that the solution of the problem of the minimization of the
functional $J_{\lambda ,\gamma }\left( u,p\right) $ in (\ref{4.5}) satisfies
the zero Neumann boundary conditions in (\ref{2.2}) as well as the initial
and terminal conditions in (\ref{2.3}) and (\ref{2.5}), we adopt the
Matlab's built-in optimization toolbox \textbf{fmincon} to minimize the
discretized form of the functional $J_{\lambda ,\gamma }\left( u,p\right) $.
The minimization was done with respect to the values of functions $u,p$ at
the grid points. The iterations of \textbf{fmincon} stop when the condition $%
|\nabla J_{\lambda ,\gamma }\left( u,p\right) |<10^{-2}$ is met. To
implement, the starting point for the iterations of \textbf{fmincon} was 
\begin{equation}
u^{(0)}(x,t)=u_{T}(x), \quad p^{(0)}(x,t)=p_{0}(x)(1-\frac{t}{T})+p_{T}(x)%
\frac{t}{T},\quad t\in \lbrack 0,T],\quad T=1.  \label{5.4}
\end{equation}

We introduce the random noise in the initial and terminal conditions in %
\eqref{2.3} and \eqref{2.5} as follows: 
\begin{equation}
\begin{split}
& \hspace{2cm}u_{T,\zeta }(x)=u_{T}(x)\left( 1+\delta \zeta _{u,x}\right) ,
\\
& p_{0,\zeta }(x)=p_{0}(x)\left( 1+\delta \zeta _{p,0,x}\right) ,\quad
p_{T,\zeta }(x)=p_{T}(x)\left( 1+\delta \zeta _{p,T,x}\right) ,
\end{split}
\label{5.41}
\end{equation}%
where $\zeta _{u,x},\zeta _{p,0,x},\zeta _{p,T,x}$ are the uniformly
distributed random variables in the interval $[0,1]$ depending on the point $%
x\in \Omega $ with $\delta =0.03$, which corresponds respectively to the $%
3\% $ noise level. The reconstructions from the noisy data are denoted as $%
u_{\zeta }(x,t),p_{\zeta }(x,t)$.

Test 1 serves as a reference test for us. This means that we select optimal
values of parameters $\gamma $ and $\lambda $ in this test and use the same
values of these parameters in the remaining Tests 2 and 3.

\textbf{Test 1.} In this test, we generate the data by the method of
subsection 5.1 for the case when the function $v\left( x,t\right) $ is a
polynomial function with the zero Neumann boundary condition in \eqref{2.2}, 
\begin{equation}
v(x,t)=(x_{1}^{2}(x_{1}-1)^{2}(x_{1}+1))(x_{2}^{2}(x_{2}-1)^{2}(x_{2}+2))(t^{2}+1).
\label{5.5}
\end{equation}%
Then we generate functions $u(x,t)$ and $p(x,t)$ as well as the input data $%
u_{T}(x)$, $p_{0}(x)$, $p_{T}(x)$ as indicated in subsection 5.1.

To evaluate the accuracy of the solution of inverse problem, denote the
relative errors in $L_{2}$ norm as: 
\begin{equation}
u_{\text{E}}=\frac{\Vert u(x,t)-u_{\text{comp}}(x,t)\Vert _{L_{2}(Q_{T})}}{%
\Vert u(x,t)\Vert _{L_{2}(Q_{T})}},\quad p_{\text{E}}=\frac{\Vert p(x,t)-p_{%
\text{comp}}(x,t)\Vert _{L_{2}(Q_{T})}}{\Vert p(x,t)\Vert _{L_{2}(Q_{T})}}.
\label{5.6}
\end{equation}%
Here, $\left( u_{\text{comp}},p_{\text{comp}}\right) (x,t)$ is the solution
of the Minimization Problem. Since we work with the finite differences in
our computations, then $L_{2}(Q_{T})-$norms in (\ref{5.6}) are understood in
the discrete sense.

To choose an optimal value of the parameter $\lambda $ in (\ref{4.5}), the
relative errors $u_{\text{E}}$ and $p_{\text{E}}$ of (\ref{5.6}) for Test 1
are displayed in Table \ref{table_relative_error_vs_lambda}. One can observe
that 
\begin{equation}
\lambda =2  \label{50}
\end{equation}
is the optimal choice with smallest error. The reconstruction of $p$ is more
accurate than $u$, because we know in (\ref{2.3}), (\ref{2.5}) both
functions $p\left( x,0\right) $ and $p\left( x,T\right) ,$ whereas we know
only $u\left( x,T\right) $. The cross-sections with $%
x_{2}=0.2,0.5,0.8,t=0.2,0.5,0.8$ of superimposed functions $u(x,t),u_{\text{%
comp}}(x,t)$ and $p(x,t),p_{\text{comp}}(x,t)$ are displayed in Figure \ref%
{plot_poly_cross_v_x} and Figure \ref{plot_poly_cross_m_x}. The
reconstructions of functions $u$ and $p$ are accurate.

\begin{table}[htbp]
\caption{The relative errors $u_{\text{E}}, p_{\text{E}}$ depending on $%
\protect\lambda $ for Test 1.}
\label{table_relative_error_vs_lambda}\centering
\vskip - 0.3 cm 
\begin{tabular}{c|c|c|c|c|c|c|c}
\hline
$\lambda$ & 0.01 & 0.5 & 1 & 2 & 3 & 4 & 6 \\ \hline
$u_{\text{E}}$ & 0.2275 & 0.2015 & 0.1536 & 0.1153 & 0.4738 & 0.5969 & 0.6070
\\ \hline
$p_{\text{E}}$ & 0.0556 & 0.0513 & 0.0482 & 0.0303 & 0.1068 & 0.1290 & 0.1301
\\ \hline
\end{tabular}%
\end{table}

\begin{figure}[htbp]
\centering
\includegraphics[width = 5in]{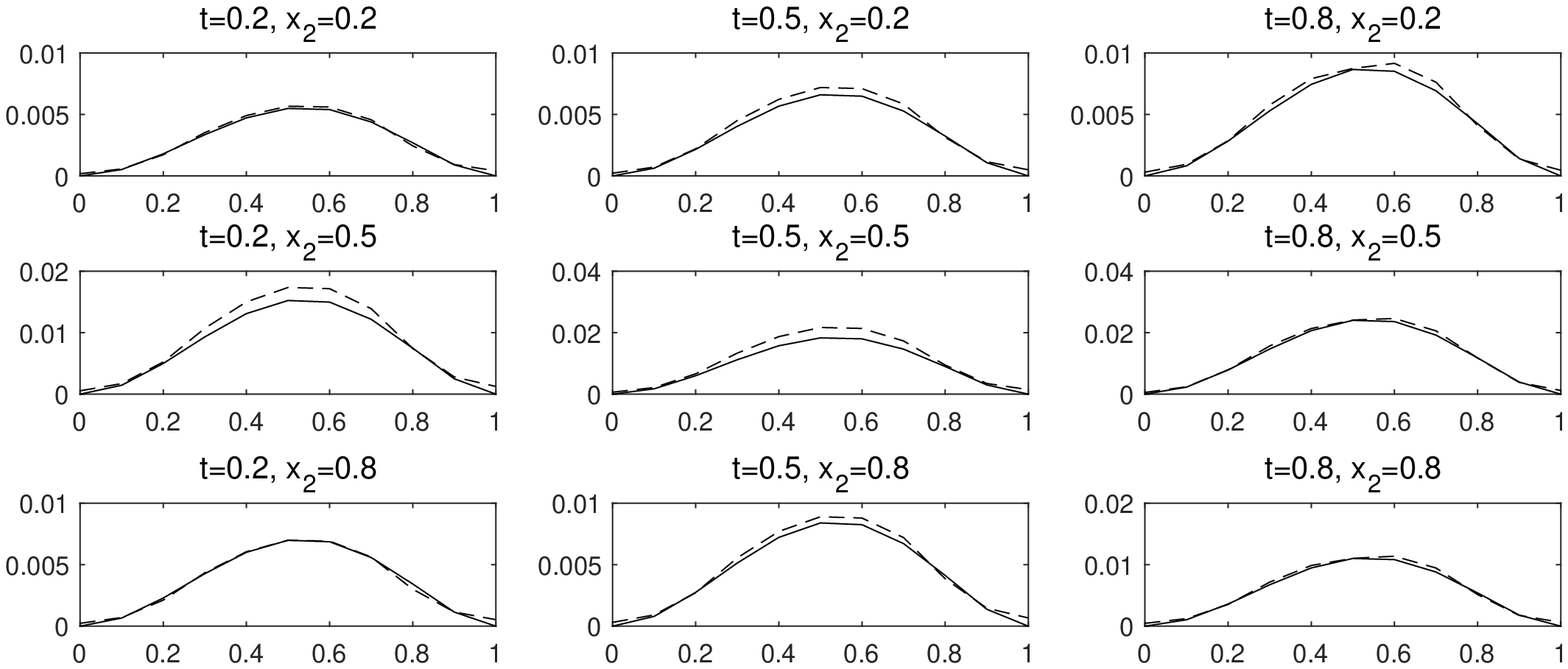} \vskip - 0.3
cm
\caption{The cross-sections of $u(x, t)$ (solid line) generated by the
procedure of subsection 5.1 with the polynomial function $v(x, t)$ in (%
\protect\ref{5.5}), and corresponding reconstruction $u_{\text{comp}}\left(
x,t\right) $ (dotted line) on the set $\{x_{2}, t=0.2, 0.5, 0.8 \} $ with $%
\protect\lambda=2$ in (\protect\ref{3.1}).}
\label{plot_poly_cross_v_x}
\end{figure}

\begin{figure}[htbp]
\centering
\includegraphics[width = 5in]{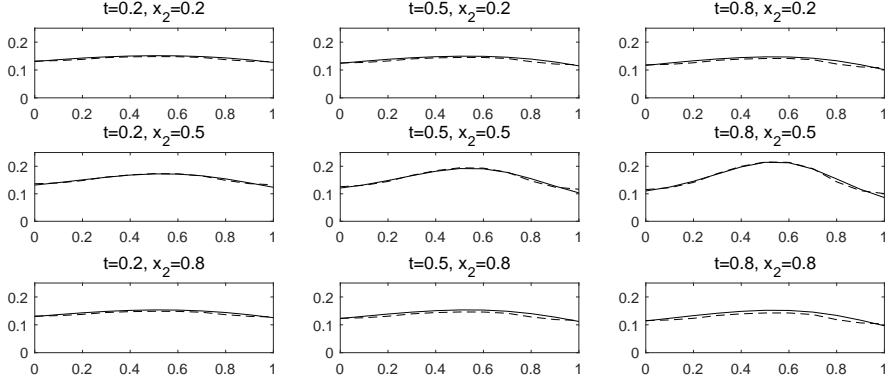} \vskip - 0.3
cm
\caption{The cross-sections of $p(x, t)$ (solid line) generated by the
procedure of subsection 5.1 with the polynomial function $v(x, t)$ in (%
\protect\ref{5.5}), and corresponding reconstruction $p_{\text{comp}}\left(
x,t\right)$ (dotted line) on the set $\{ x_{2}, t=0.2, 0.5, 0.8 \} $ with $%
\protect\lambda=2$ in (\protect\ref{3.1}).}
\label{plot_poly_cross_m_x}
\end{figure}

\textbf{Test 2.} In this test, we use a trigonometric function $v(x,t)$ with
the zero Neumann boundary condition for the data generation procedure of
subsection 5.1, 
\begin{equation}
v(x,t)=0.01\cos (\pi x_{1})\sin (\pi (x_{2}-0.5))(t^{2}+1).  \label{5.7}
\end{equation}%
The evolutions of functions $u(x,t)$ and $p(x,t)$ at $t=0,0.2,0.5,0.8,1$,
generated by the procedure of subsection 5.1 is exhibited in Figure \ref%
{plot_sin_v_m_t}. The cross-sections with $x_{2}=0.2,0.5,0.8,t=0.2,0.5,0.8$
of $u(x,t),u_{\text{comp}}(x,t)$ and $p(x,t),p_{\text{comp}}(x,t)$ are shown
in Figure \ref{plot_sin_cross_v_x} and Figure \ref{plot_sin_cross_m_x} with $%
\lambda =2$ as in (\ref{50}). The reconstructions of $u$ and $p$ are
accurate.

\begin{figure}[htbp]
\centering
\includegraphics[width = 5in]{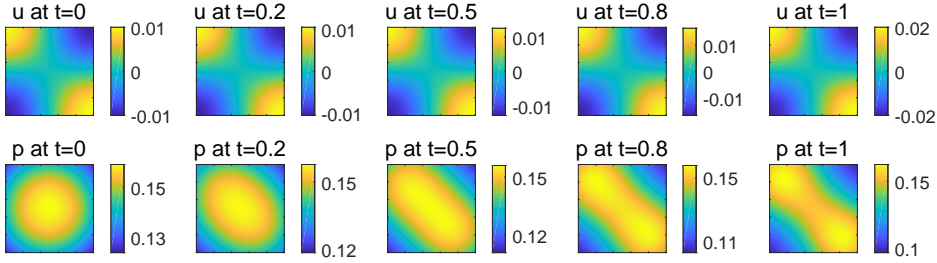} \vskip - 0.3 cm
\caption{Evolutions of functions $u\left( x,t\right) $ and $p(x, t)$
generated by the procedure of subsection 5.1 with the trigonometric function 
$v(x, t)$ in (\protect\ref{5.7}).}
\label{plot_sin_v_m_t}
\end{figure}

\begin{figure}[htbp]
\centering
\includegraphics[width = 5in]{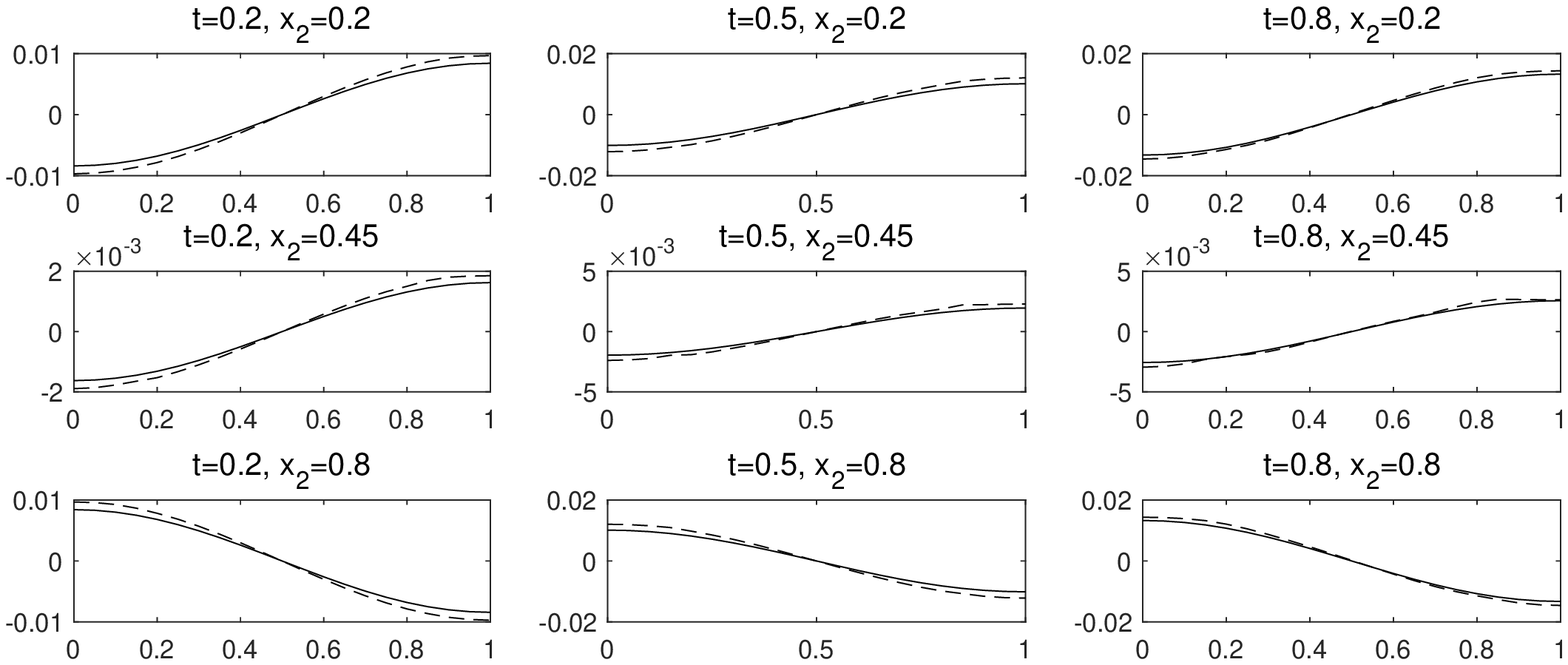} \vskip - 0.3 cm
\caption{The cross-sections of $u(x, t)$ (solid line) generated by the
procedure of subsection 5.1 with the trigonometric function $v(x, t)$ in (%
\protect\ref{5.7}), and corresponding reconstruction $u_{\text{comp}}(x, t)$
(dotted line) on the set $\left( x_{2},t\right) =\left\{
0.2,0.45,0.8\right\} \times \left\{ 0.2,0.5,0.8\right\}$.}
\label{plot_sin_cross_v_x}
\end{figure}

\begin{figure}[htbp]
\centering
\includegraphics[width = 5in]{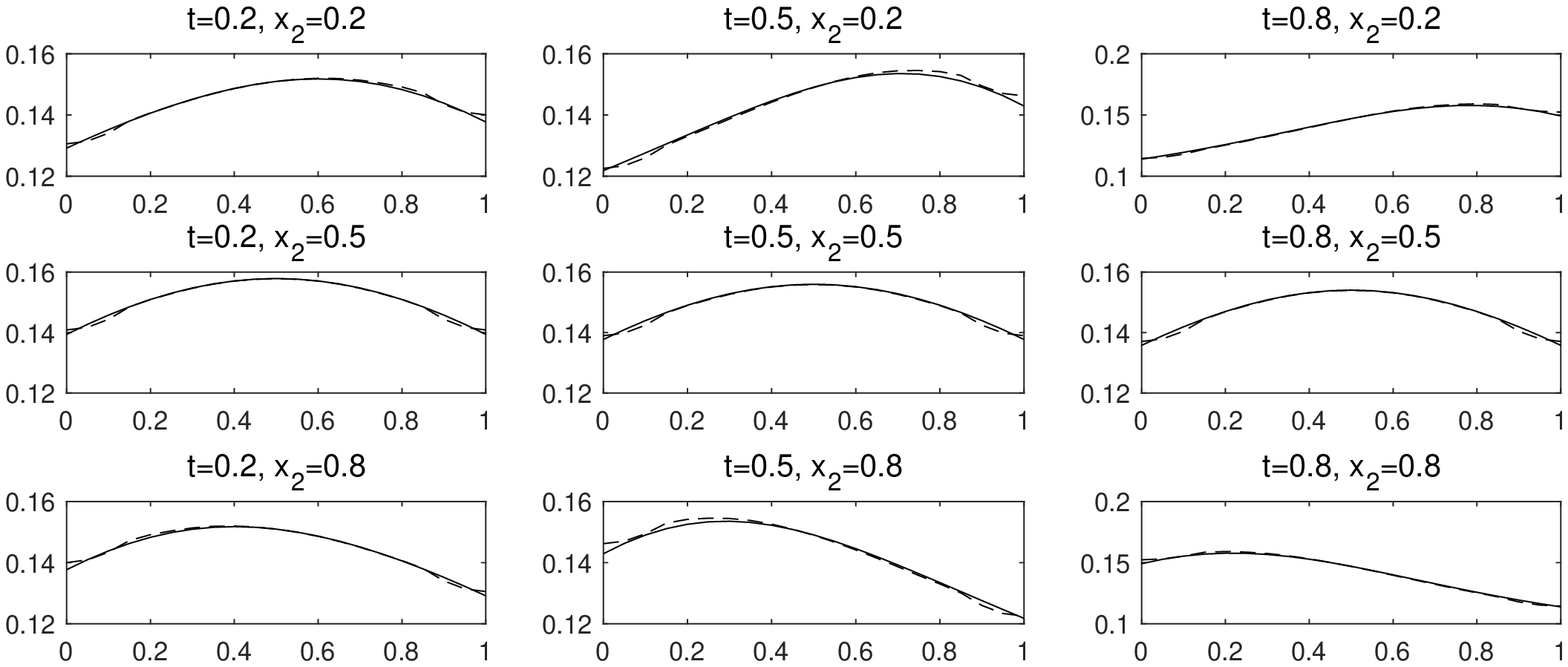} \vskip - 0.3 cm
\caption{The cross-sections of $p(x, t)$ (solid line) generated by the
procedure of subsection 5.1 with the trigonometric function $v(x, t)$ in (%
\protect\ref{5.7}), and corresponding reconstruction $p_{\text{comp}}(x, t)$
(red line) on the set $\left( x_{2},t\right) =\left\{ 0.2,0.45,0.8\right\}
\times \left\{ 0.2,0.5,0.8\right\}$.}
\label{plot_sin_cross_m_x}
\end{figure}

\textbf{Test 3.} We now test the ability of our method to work with rather
complicated non-convex shapes for both noiseless and noisy data. We again
use $\lambda =2$ as in (\ref{50}). In this test, the case when the function $%
v(x,t)$, which generates our noiseless data as in subsection 5.1, has the
shape of the letter `$C$' for each $t,$ and the size of this letter expands
when $t$ increases from $0$ to $1$. The function $v\left( x,t\right) $ is
given by:%
\begin{equation}
\begin{split}
& \hspace{0.5cm}v(x,t)=f_{v}\left( d(x),t\right) (1+t^{2}),\quad d(x)=\sqrt{%
(x_{1}-0.6)^{2}+(x_{2}-0.5)^{2}}, \\
& \hspace{1cm}f_{v}\left( d(x),t\right) =\left\{ 
\begin{array}{cc}
h_{v}\left( d(x),t\right) , & r_{1}\leq d\leq r_{2},\text{ }x_{1}\leq 0.75,\ 
\\ 
0, & \text{otherwise},%
\end{array}%
\right. \\
& \hspace{0.01cm}h_{v}\left( d(x),t\right)
=0.1(0.75-x_{1})(d(x)-r_{1})(r_{2}-d(x))\exp \left(
-100(d(x)-r_{3})^{2}\right) , \\
& r_{1}(t)=0.05(1-t)+0.15t,\quad r_{2}(t)=0.35(1-t)+0.45t,\quad
r_{3}=(r_{1}+r_{2})/2.
\end{split}
\label{5.8}
\end{equation}

We display the results with $t=0,0.2,0.5,0.8,1$ in Figure \ref%
{plot_C_moving_re_v} and Figure \ref{plot_C_moving_re_m}. First rows are the
exact functions $u\left( x,t\right) $ and $p\left( x,t\right) $ generated by
the procedure of subsection 5.1, second rows are the reconstructed function $%
u_{\text{comp}}\left( x,t\right) $ and $p_{\text{comp}}\left( x,t\right) $
with noiseless data, third rows are the reconstructed function $u_{\zeta
}\left( x,t\right) $ and $p_{\zeta }\left( x,t\right) $ for the case of
noisy data with 3\% noise level in (\ref{5.41}). The reconstructions of
functions $u\left( x,t\right) $ and $p\left( x,t\right) $ for noiseless and
noisy data are all accurate.

\begin{figure}[htbp]
\centering
\includegraphics[width = 5in]{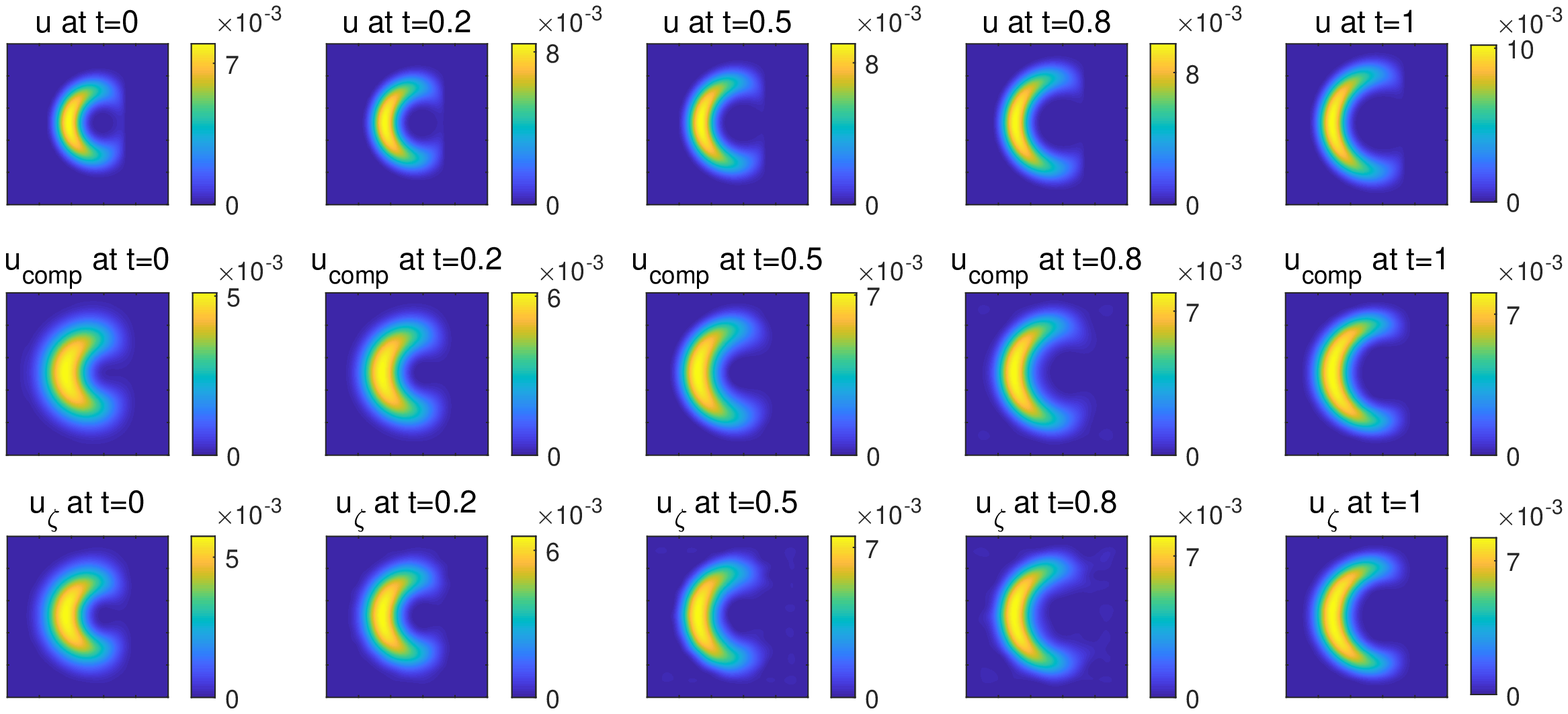} \vskip - 0.3 cm
\caption{Evolutions of generated and computed functions$u\left( x,t\right) $
for Test 3 with the function $v\left(x,t\right) $ having the shape of the
letter `C', the size of which is expanding as $t$ is increasing from 0 to 1.
The function $v\left(x,t\right) $ is given by formula (\protect\ref{5.8}).
First raw: exact function $u\left( x,t\right) $ generated by the procedure
of subsection 5.1 for $t=0,0.2,0.5,0.8,1$. Second raw: the reconstructed
function $u_{\text{comp}}\left( x,t\right) $ for case of noiseless data.
Third raw: the reconstructed function $u_{\protect\zeta }\left( x,t\right)$
for the case of noisy data with 3\% noise level. The noisy data are
generated by (\protect\ref{5.41}). One can observe accurate reconstructions
of the function $u\left( x,t\right) $ for both noiseless and noisy cases.}
\label{plot_C_moving_re_v}
\end{figure}

\begin{figure}[tbph]
\centering
\includegraphics[width = 5in]{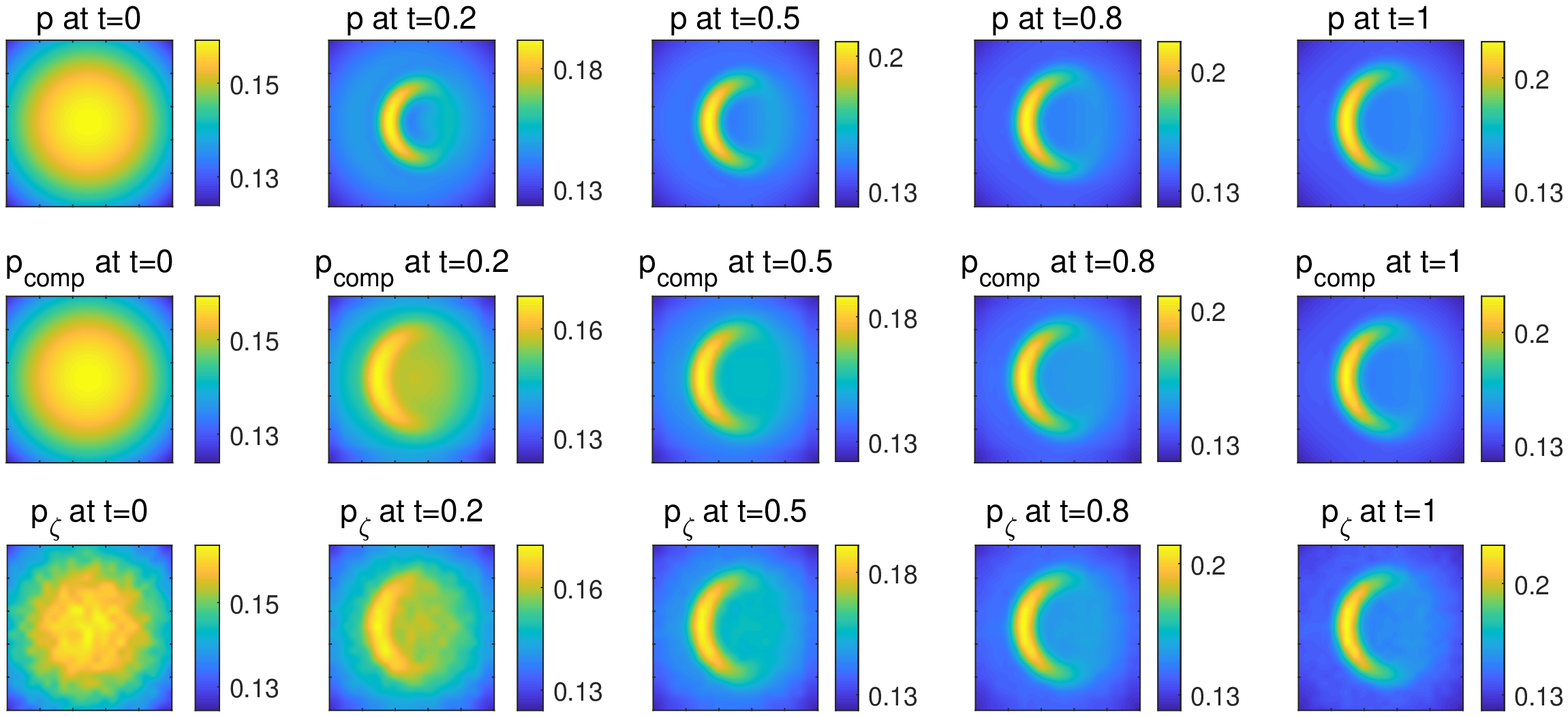} \vskip - 0.3 cm
\caption{Evolutions of generated and computed functions $p\left( x,t\right) $
for Test 3. First raw: exact function $p\left( x,t\right) $ generated by the
procedure of subsection 5.1 for $t=0,0.2,0.5,0.8,1$. Second raw: the
reconstructed function $p_{\text{comp}}\left( x,t\right) $ for the case of
noiseless data. Third raw: the reconstructed function $p_{\protect\zeta %
}\left( x,t\right) $ for the case of noisy data with 3\% noise level. The
noisy data are generated by (\protect\ref{5.41}). One can observe accurate
reconstructions of the function $p\left( x,t\right) $ for both noiseless and
noisy data.}
\label{plot_C_moving_re_m}
\end{figure}

\section{Summary}

\label{sec:6}

We have developed in this paper the first numerical method for the Mean
Field Games System of the second order with rigorously guaranteed global
convergence. We have addressed the problem of the retrospective analysis for
this system. Convergence rates are explicitly given, and their depend on the
iteration number and the level of noise in the data. The key idea is the
idea of applications of the tool of Carleman estimates and, based on it, the
convexification numerical method. In the previous works of this research
group, this method has consistently demonstrated a good performance for a
variety of nonlinear ill-posed problems and coefficient inverse problems 
\cite{Bak,SAR,Kpar,KL,Ktransp,LeLoc1,LeLoc2}.

The idea of using Carleman estimates was first introduced in the MFG theory
in the work \cite{MFG1} and continued since then in a number of publications
of this research group, see \cite{MFG2,MFG6} and references cited in \cite%
{MFG6}. Results of numerical studies demonstrate a good accuracy of our
reconstructions for both noiseless and noisy data.

\section*{Acknowledgement}

We thank Professor Jie Xiong from SUSTech for many fruitful discussions which lead to significant improvement of the manuscript.


\vskip 1 cm

\end{document}